\documentclass{arxiv}[11pt]
\pdfoutput=1

\usepackage{amssymb}

\usepackage{times}
\usepackage{ae, aecompl} 
\usepackage{textcomp}
\usepackage[english]{babel}
\usepackage{natbib}

\usepackage{etex}

\usepackage{amsmath} 
\usepackage{amsthm} 
\usepackage{amsfonts} 
\usepackage{mathdots} 
\usepackage{url} 
\usepackage{ltxtable} 
\usepackage{graphicx} 
\usepackage[table]{xcolor} 
\usepackage{eurosym} 
\usepackage[active]{srcltx}
\usepackage{pifont}
\usepackage{snapshot}
\usepackage{pseudocode}

\usepackage{animate}
\usepackage{multirow}
\usepackage{bigstrut}
\usepackage{xcolor}
\usepackage{booktabs,colortbl,tabularx}

\usepackage{array}
\usepackage{subfig}

\newcommand{\ra}[1]{\renewcommand{\arraystretch}{#1}}

\newcommand{\rod}{\mathbf{D}}

\newcommand{\es}{\pmb{\tau}}
\newcommand{\vel}{\mathbf{u}}

\newcommand{\divergence}{\nabla \cdot}

\newcommand{\cauchystress}{\pmb{\sigma}}

\newcommand{\Rey}{\text{Re}\,}

\newcommand{\spectralhp}{spectral/$hp$}

\newcommand{\mytoprule}{\specialrule{1.5pt}{0em}{0em}}
\newcommand{\mymidrule}{\specialrule{1pt}{0em}{0em}}
\newcommand{\mybottomrule}{\specialrule{1.5pt}{0em}{0em}}

\newtheorem{theorem}{Theorem}[section] 
\newtheorem{Problem}[theorem]{Problem}

\title{Spectral/hp element methods for plane {N}ewtonian extrudate swell}

\title{Spectral/hp element methods for plane {N}ewtonian extrudate swell}
 \author{S. Claus$^1$, C. D. Cantwell$^2$ and
  T.N. Phillips$^3$} 

\address{$^{1}$ Department of Mathematics, University College London, Gower Street,
London WC1E 6BT, United Kingdom, susanne.claus@ucl.ac.uk\\
  $^{2}$ Department of Aeronautics, Roderic Hill Building, Imperial College, London SW7 2AZ, United Kingdom \\
  $^{3}$ School of Mathematics, Cardiff University, Senghennydd Road, Cardiff CF24 4AG, United Kingdom}

\begin{document}

\begin{summary}
Spectral/hp element methods and an
arbitrary Lagrangian-Eulerian (ALE) moving-boundary technique are used to investigate planar Newtonian extrudate swell. Newtonian extrudate swell arises when viscous liquids exit long die slits. The problem is 
characterised by a stress singularity at the end of the slit which is
inherently difficult to capture and strongly influences the predicted
swelling of the fluid. The impact of inertia ($0 \leq \Rey \leq 100$) and slip along the die wall on the free surface profile and the velocity and pressure values in the domain and around the singularity are investigated. The high order method is shown to provide high resolution of the steep pressure profile at the singularity. The swelling ratio and exit pressure loss are compared with existing results in the literature and the ability of high-order methods to capture these values using significantly fewer degrees of freedom is demonstrated.
\end{summary}




\section{Introduction}
\label{sec: Introduction}
In this article, we investigate the extrudate swell phenomenon, which is a
radial swelling of free liquid jets exhibited by viscous fluids exiting long die
slits. This jet swelling is particularly strong for viscoelastic fluids but is also
exhibited by low Reynolds number Newtonian fluids. The prediction of the
swelling ratio is very important in a range of industrial processes such as
inkjet printing, extrusion moulding or cable coating.

The swelling of Newtonian jets is mainly characterised by the reorganisation of the
velocity profile from the parabolic Poiseuille flow inside the die to plug flow
downstream (\cite{Tanner:2002}). This transition is characterised by the sudden
jump in the shear stress at the die exit (\cite{Russo2009}). Inside the die, the
shear stress at the wall is at its maximum with particles sticking to the wall
(for the no-slip boundary condition). Then immediately after the die exit, the
removal of the wall shear stress causes a boundary layer to form at the free
surface. In this layer, the parabolic velocity profile adjusts itself so as to
satisfy the condition of zero shear stress at the free surface. This sudden jump
in the shear stress at the die exit causes an almost instantaneous acceleration
of the particles at the free surface causing the fluid jet to swell.

Due to the presence of this stress singularity at the die exit, numerical
simulations of the extrudate swell phenomenon are particularly challenging.
Analytically, this singularity originates from the sudden change in the boundary
condition from the wall of the die to the free surface of the exiting jet. This
"jump" in the boundary condition yields steep and infinite stress and pressure
concentrations at the singular point. These infinite stress values near the
singularity affect the accuracy of the numerical solution and the size of the
swelling and therefore need to be resolved as accurately as possible. In this
contribution, we use a \spectralhp{} element method to improve our ability to
capture these stress concentrations. Traditional discretisation methods such as
finite differences or low-order finite elements require a very large number of
degrees of freedom to resolve these sharp stress variations.

In this article, we will describe a \spectralhp{} method that is capable of
approximating the infinite stress values with an exponential increase in the
extreme values of the pressure with $p$-refinement. This demonstrates that our
high order method provides a high-quality approximation of the stress
singularity with a very low number of degrees of freedom. We will give detailed
information about the pressure and velocity in the vicinity of the singularity
for a wide range of Reynolds numbers ($0 \leq \Rey \leq 100$) and for slip along
the die wall. We demonstrate that our method predicts swell ratios and exit
pressure loss corrections in excellent agreement with a recent numerical study
of \cite{Mitsoulis2012a} for our coarsest approximation $P=10$.
\cite{Mitsoulis2012a}  used a low order finite element method with a high mesh
refinement around the singularity.

Typically, a decrease in the swelling is observed for an increase in the
resolution of the singularity. In the existing literature, the stress values at
the singularity are rarely addressed.   \cite{Salamon1995} investigated the role of surface tension and slip on the singularity  numerically and analytically. They demonstrated that a very fine mesh near the singularity is needed to predict the singular pressure and stress behaviour with sufficient accuracy.
\cite{Georgiou1999}  compared the
singular finite element method with the regular finite element method for the extrudate swell problem. In the singular finite element method basis functions in the elements
around the singularity are enriched with the local asymptotic solution for the
singularity. They demonstrated that with this method the
predictions of the swell ratio converged. However, the singular finite element
method requires the correct asymptotic behaviour of the pressure at the corner
singularity and the asymptotic solution for the pressure is obtained assuming
Stokes-like behaviour around the singularity. This means this approach is only
accurate for $\Rey=0$. Indeed, \cite{Georgiou1999} found that the singular
finite element method was outperformed by the regular finite element method for
extrudate swell including inertia. Our method is capable of resolving the stress
singularity with spectral convergence properties without making any assumptions
on the form of the singularity.

Inertialess extrudate planar Newtonian swell has been investigated in terms of
swell ratios using low order finite elements by a wide range of authors 
(\cite{Tanner1973}, \cite{Nickell1974}, \cite{Crochet1982}, \cite{Reddy1978}).
\cite{Tanner:2002} provides a review of inertialess Newtonian swell ratio
results.  Only very few investigations involved the use of  higher order
methods. \cite{Ho1994} provided the first extrudate swell computation with a
spectral method for one coarse mesh with $8$ spectral elements with polynomial
order $4$ for $\Rey=0$. They predicted a swell ratio of $1.1840$.
\cite{Russo2009} used the spectral element method to predict free surface
profiles and swell ratios for $0\leq \Rey \leq 10$ and surface tension for $4$
spectral elements with polynomial order $6<P<14$. We will use a spectral element
mesh with $14$ spectral elements and $10 \leq P \leq 16$ with a smaller element
size around the singularity providing a much higher resolution there compared
with previous studies. We will provide results for $0\leq \Rey \leq 100$ and for
a slip condition along the die wall.

The paper is organised as follows. In Section~\ref{sec: Formulation}, we will introduce the governing equations for the description of Newtonian free surface flow and the equations of motion for the mesh movement. We will conclude this Section with a description of the boundary conditions for the extrudate swell problem and the definition of the quantities of interest such as swelling ratio and exit pressure correction. In Section~\ref{sec: Numerical Discretisation}, we describe the numerical discretisation of the governing equations. In Section~\ref{sec: numerical results}, we give numerical results for the impact of inertia and slip on the extrudate swell problem including detailed plots for velocity and pressure profiles in different parts of the domain.

\section{Formulation}
\label{sec: Formulation}
\subsection{Governing Equations of the Fluid}
\label{sec: Governing Equations}
The free surface motion of an incompressible fluid flow can be characterised by
the incompressible Navier-Stokes equations describing the motion of the fluid
and the motion of the free surface.  On a moving domain 
$\Omega_t \subset \mathbb{R}^d$, $t \in I \equiv (t_0,T)$, they can be expressed
as
\begin{subequations}
\label{equ:navier-stokes}
\begin{align}
\Rey \left(\dfrac{\partial \mathbf{u}}{ \partial t} 
            + (\mathbf{u} \cdot \nabla)\mathbf{u}\right)
    &= -\nabla p + 2 \, \divergence \rod  
    &&\mbox{in } \Omega_t, t\in I, \\
\divergence \mathbf{u} &= 0  
    &&\mbox{in } \Omega_t, t\in I, \\ 
\vel &= \vel_0  &&\mbox{in } \hat{\Omega}_{t_0}, \\
\vel &= \vel_D  &&\mbox{on } \partial \Omega_t, t  \in I,
\end{align}
\end{subequations}
where $\mathbf{u}$ is the velocity, $p$ is the pressure,  
$\mathbf{D} = \frac12 ( \nabla \vel +  \nabla \vel^T)$ is the rate of deformation 
tensor, $\Rey$ is the Reynolds number, $\vel_0$ is the  velocity field at $t=t_0$
and $\vel_D$ is the assigned Dirichlet boundary condition.

The motion of the free surface, $\Gamma_f$, is characterised by the following
boundary conditions
\begin{subequations}
\begin{align}
\vel \cdot \mathbf{n} &= \mathbf{w} \cdot \mathbf{n}  
    &&\mbox{on } \Gamma_f \quad \mbox{ (kinematic)} \label{equ:kinematic-bc}\\
 \left[\cauchystress\right] \cdot \mathbf{n} &= \sigma \kappa \mathbf{n}
    &&\mbox{on } \Gamma_f \quad \mbox{ (dynamic)} \label{equ:dynamic-bc}
\end{align}
\end{subequations}
where $\mathbf{w}$ is the velocity of the free surface, $\sigma$ is the surface
tension coefficient, $\kappa$ is the curvature of the free surface, $\mathbf{n}$
is the unit outward normal to the free surface and  $\left[\cauchystress\right]$
denotes the jump in the Cauchy stress tensor across the free surface.

In order to track the free surface motion computationally, the grid
points of our computational mesh at the free surface are moved with the normal fluid
velocity, which ensures that particles do not cross the interface and therefore that the kinematic condition~\eqref{equ:kinematic-bc}
is satisfied.

To avoid mesh distortion, the mesh points in the interior of the domain
are moved with an arbitrary speed. This use of arbitrary mesh movement is  known as the arbitrary
Lagrangian-Eulerian (ALE) technique. The ALE formulation relates the
Navier-Stokes equations on the moving domain~\eqref{equ:navier-stokes} to a
formulation on a referential configuration $\hat{\Omega}_{t_0}$. At each $t\in
I$, each point of the reference configuration $\mathbf{Y}$ is then associated to
a point $\pmb{x}$ in the current domain $\Omega_t$ using the so-called ALE-map \citep{Donea2004, Scovazzi2007, Pena2009,Nobile2001}, that is,
\begin{align}
\mathcal{R}_{t}: \hat{\Omega}_{t_0} \rightarrow &\Omega_{t}, 
    && \forall t \geq 0 , \nonumber \\
\mathbf{Y} \mapsto &\pmb{x}(\mathbf{Y},t)=\mathcal{R}_{t}(\mathbf{Y}), 
    && \forall \mathbf{\mathbf{Y}} \in \hat{\Omega}_{t_0},
\label{equ: referential maprepeat}
\end{align}  
where $\mathbf{Y}$ is called the ALE coordinate and $\pmb{x}$ is the Eulerian 
coordinate. The movement of the mesh, can then be characterised by the 
following quantities
\begin{enumerate}
\item the mesh velocity 
\begin{align}
\mathbf{w}(\pmb{x},t)
:= \left.\dfrac{\partial \pmb{x}(\mathbf{Y},t)}{\partial t}\right|_{\mathbf{Y}} 
= \left.\dfrac{\partial \mathcal{R}_{t}(\mathbf{Y}) }{\partial t} \right|_{\mathbf{Y}}.
\label{equ: mesh velocityrepeat}
\end{align}
\item the material time derivative in terms of the time derivative with respect to the ALE-frame
\begin{align}
\dfrac{Df(\pmb{x},t)}{Dt}=\left. \dfrac{\partial f}{\partial t} \right|_{\mathbf{Y}}  + \left(\vel-\mathbf{w}\right) \cdot \nabla_{\pmb{x}} f 
\label{equ: ALEtimederivativerepeat}
\end{align}
\end{enumerate}
Equation~\eqref{equ:navier-stokes} in the ALE-formulation reads
\begin{subequations}
\label{equ: conservation equations ALE-form}
\begin{align}
\Rey \left(\left. \dfrac{\partial \vel}{\partial t} \right|_{\mathbf{Y}} 
            + \left(\vel - \mathbf{w}\right) \cdot \nabla_{\pmb{x}} \vel\right)
    &= -\nabla_{\pmb{x}}p + 2\,\nabla_{\pmb{x}}\cdot\rod_{\pmb{x}}
  &&\mbox{for } \pmb{x} \in \Omega_t, t  \in I, \\
\nabla_{\pmb{x}} \cdot \vel &= 0, 
    &&\mbox{for } \pmb{x} \in \Omega_t,\, t  \in I,\\
\vel &= \vel_0,  
    &&\mbox{for } \pmb{x} \in\hat{\Omega}_{t_0},\\
\vel &= \vel_D,   
    &&\mbox{on }  \pmb{x} \in\partial \Omega_t, \, t  \in I.
\end{align}
\end{subequations}
Here, $\mathbf{D}_{\pmb{x}} = \frac12 ( \nabla_{\pmb{x}} \vel +  \nabla_{\pmb{x}} 
\vel^T)$ is the rate of deformation tensor in the Eulerian frame of reference.

\subsection{Governing Equations of the Mesh}
In addition to the motion of the fluid, we need to find a sensible way to
describe the domain movement. In general, the domain movement is characterised
by the movement of its boundary $\partial \Omega_t$ and can be described using
the domain or mesh velocity $\mathbf{w}$ (\cite{Ho1994}, \cite{Robertson2004}),
the ALE-mapping $\mathcal{R}(t)$ (\cite{Nobile2001}, \cite{Pena2009}) or the
displacement $\mathbf{d}= \Delta t \mathbf{w}$ (\cite{Choi2011}). In the present
work, we describe the domain movement using the mesh velocity, $\mathbf{w}$. For
the domain movement, we choose boundary conditions such that the kinematic
boundary condition is satisfied and mesh distortions are kept to a minimum, that
is,
\begin{subequations}
\label{equ: mesh velocity boundary conditions}
\begin{align}
\mathbf{w} \cdot \mathbf{n} &= \mathbf{u} \cdot \mathbf{n},  
&& \mbox{on } \Gamma_f(t),  \\
\mathbf{w} \cdot \mathbf{s} &= 0
&& \mbox{on } \Gamma_f(t),  \\
\nabla \mathbf{w} \cdot \mathbf{n} &= \mathbf{0} 
&& \mbox{at outflow}, \\
\mathbf{w} &= \mathbf{0} 
&& \text{elsewhere},
\end{align}
\end{subequations}
where $\mathbf{s}$ is the unit tangential vector on the free surface boundary. In 
order to guarantee smooth mesh movement in the interior, we solve an elliptic 
problem for the mesh velocity, given by 
\begin{align}
\Delta \mathbf{w} = 0 \mbox{ on } \Omega(t).
\label{equ: mesh velocity strong form}
\end{align} 
subject to the boundary conditions \eqref{equ: mesh velocity boundary conditions}. This harmonic mesh movement preserves a high quality mesh for small displacements and has been employed, for instance, by \cite{Ho1994},
\cite{Nobile2001} and \cite{Pena2009}. However, for higher mesh deformations,
other elliptic problems may be solved for the movement of the domain, such as
elliptic operators arising from Stokes or elasticity problems (see the monograph
of \cite{Deville2002} for further details).

\subsection{Computational Domain and Quantities of Interest}
\label{sec: computational domain swell}

\begin{figure*}[t]
\centering
\includegraphics[width=\linewidth]{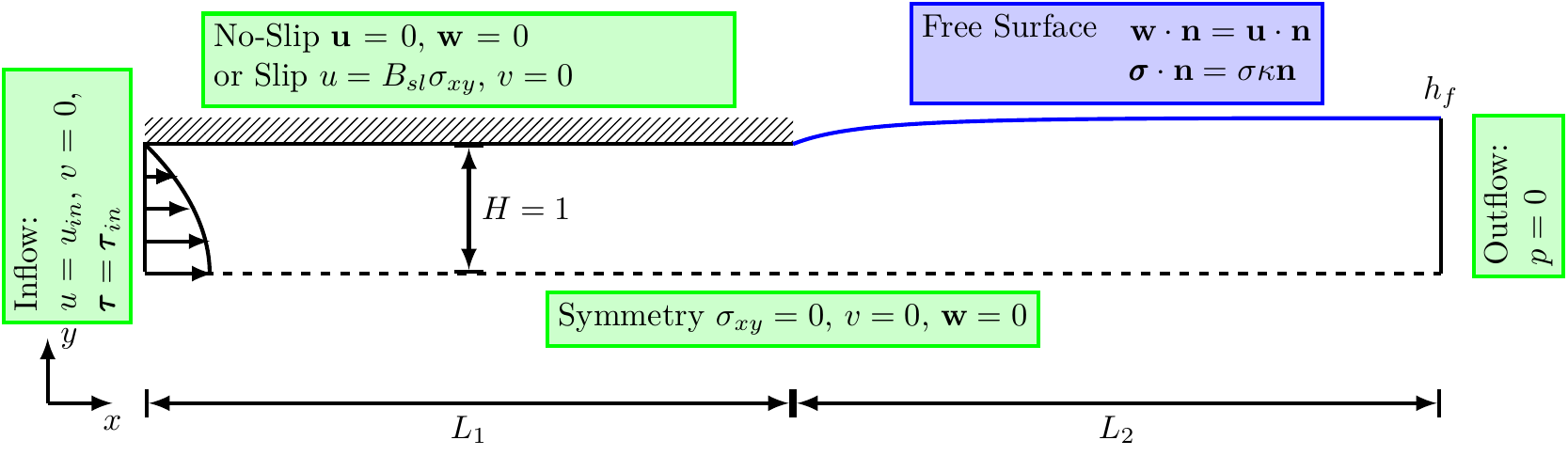}
\caption{Schematic of the symmetric die swell flow configuration. $L_1$ is the
length of the die which has fixed boundaries and is of half-height H.
$L_2$ is the length of the outflow region, the boundaries of which are
free to move. Boundary conditions are provided for each surface.}
\label{fig:DieSwell}
\end{figure*}

Consider the extrusion of a Newtonian liquid from a planar die. The schematic diagram
of the employed planar die geometry is depicted in Figure~\ref{fig:DieSwell}. We
consider a die of length $L_1$ and height $H$, and an exit region of length
$L_2$. The length of the die is chosen sufficiently long in order to guarantee a
fully developed flow far upstream of the exit plane.
In the following, we pay special attention to the following two quantities of
interest: the swelling ratio and the pressure exit correction factor. In
practice, the extrudate swell ratio is of importance in extrusion processes and
the excess pressure loss gives an indication as to how much extra pressure has to be
applied to achieve certain swell ratios. The swelling ratio, $\chi_R$,  is
defined as
\begin{align}
\chi_R = \dfrac{h_f}{H}
\end{align} 
where $H$ is the half-height of the die and $h_f$ is the half-height of the liquid jet at the outflow boundary. The swelling ratio is a 
function of several parameters
\begin{align}
\chi_R(H,\langle u\rangle,\Rey,B_{sl}),
\end{align}
where $\langle u \rangle$ is the average 
inflow velocity, $\Rey$ is the Reynolds number and $B_{sl}$ is the slip parameter
along the die wall.

The dimensionless pressure exit correction factor,  $n_{ex}$, is 
defined as 
\begin{align}
n_{ex}=\dfrac{\Delta p - \Delta p_0}{2 \sigma_w}
\label{equ: pressure exit correction def}
\end{align}
where $\Delta p$ is the pressure drop between the inlet and the outlet plane, 
$\Delta p_0$ is the pressure drop between the inlet and the exit of the die for
fully developed Poiseuille flow and $\sigma_w$ is the shear stress at the 
channel wall corresponding to fully developed Poiseuille flow. Here, the 
pressure differences are taken along the centreline. In particular, the 
pressure drops are given by (\cite{Tanner:2002})
\begin{align}
\Delta p_0 =& p|_{x=-L_1} = 2 \sigma_w \frac{L_1}{H} 
    &&\text{Poiseuille flow for } x \in [-L_1, \, 0]\\
\Delta p =& p|_{x=-L_1} - p|_{x=L_2} 
    &&\text{Extrudate Swell for } x \in [-L_1, \, L_2].
\end{align} 

In our computations, we employ the following boundary conditions as depicted in
Figure~\ref{fig:DieSwell} for a half-channel height of $H=1$. We assume the flow
is symmetric and along the symmetry line, we set $v=0$ and $\sigma_{xy}=0$.
Note that, $\sigma_{xy}=0$ is set through the boundary integral in the momentum
equation~\eqref{equ: ALE momentum}. For the die swell geometry this means that
there is no contribution of the Neumann boundary integral in the momentum
equation along the symmetry line. At the die wall, we either impose no-slip
boundary conditions (i.e. $\vel=0$) or Navier's slip condition. The latter is a
mixed boundary condition of Dirichlet and Neumann type. For the extrudate swell
geometry depicted in Figure~\ref{fig:DieSwell}, we set $v=0$ and impose
$\sigma_{xy}=\frac{1}{B_{sl}}u$ through the Neumann boundary term in the
momentum equation~\eqref{equ: ALE momentum}. This means for the velocity
component $u$ along the slip boundary $\Gamma_{sl}$, we obtain the boundary
integral
\begin{align}
\int_{\Gamma_{sl}} \left(\cauchystress \cdot \mathbf{n} \, \phi_{\vel} \right)
    \mathbf{e}_{x}\, \mathrm{d} \Gamma =  \int_{\Gamma_{sl}} \frac{1}{B_{sl}}u 
    \phi_u \, \mathrm{d}\Gamma
\label{equ: slip boundary}
\end{align}
where $\mathbf{e}_x$ is the unit vector in $x$-direction. At the outflow, we 
employ an open outflow boundary condition. We assume a reference pressure
$p=0$ along the outflow boundary and the remaining terms in the Neumann 
boundary integral along the outflow boundary in the momentum equation are 
evaluated along with the volume integrals. 
At inflow, we either impose the parabolic profile 
\begin{align}
u=\dfrac{3}{2}\left(1- y^2\right), \quad v=0
\label{equ: inflow swell}
\end{align}
in combination with no-slip along the die wall or the profile  (\cite{Kountouriotisa:2012})
\begin{align}
u = \dfrac{3}{2(1+3B_{sl})}(1-y^2+2B_{sl}), \quad \dfrac{\partial u}{\partial y}=\dfrac{-3y}{(1+3B_{sl})}, \quad v=0
\label{equ: slip inflow}
\end{align}
in combination with the slip boundary condition. In the extrudate swell problem the velocity field undergoes a transition from Poiseuille flow inside the die to plug flow in the free jet. 
Due to the conservation of energy the flow
rate in the die has to be the same as in the uniform plug flow, which yields
\begin{align}
u_{\mathrm{plug}} = \dfrac{1}{2 h_{\mathrm{plug}}} \int_{-H}^H u(y) \,
\mathrm{d}y
\label{equ: uplug}
\end{align}
where $h_{\mathrm{plug}}$ is the height of the fluid jet in the uniform flow
region and $u(y)$ is the parabolic Poiseuille flow profile. We have 
$0<u_{\mathrm{plug}}<u_{\mathrm{max}}$, which means that while particles at the
free surface accelerate when exiting the die the flow near the centreline decelerates.

For the mesh velocity, we employ the following boundary conditions.
We consider the mesh to be fixed at inflow, the die wall and along the symmetry line, i.e. $\mathbf{w}=(w_x,w_y)=0$.  At the outflow boundary,
we allow the mesh to move in the $y$-direction, i.e.
$\nabla w_y \cdot \mathbf{n}=0$, and fix it in the $x$-direction, $w_x=0$. At
the free surface, we enforce the kinematic boundary condition through the mesh velocity in terms of a Dirichlet boundary condition for the mesh-velocity, i.e.
\begin{align}
\mathbf{w} \cdot \mathbf{n} = \mathbf{u} \cdot \mathbf{n}.
\end{align} 
To avoid mesh distortion, we choose to move the mesh along the free surface
boundary only in the $y$-direction. The mesh is moved with sufficient velocity
$w_y$ into the $y$-direction to ensure that no particle crosses the interface,
that is,
\begin{align}
w_x = 0 , \quad w_y =v + u \frac{n_x}{n_y}.
\end{align}

\section{Numerical Discretisation}
\label{sec: Numerical Discretisation}

\subsection{Spectral Element Discretisation}
Consider the decomposition of the domain $\Omega_t^{\delta}$ into $N_{el}$
non-overlapping elements. These elements are each mapped to
a standard element on which the unknowns are approximated using a modal
polynomial expansion basis proposed by \cite{Dubiner:1991} and extended by
\cite{Sherwin:2005} given by
\begin{align}
\phi_p(\xi) = 
\begin{cases}
  \dfrac{1-\xi}{2}, \quad & p=0,\\[1.5ex]
  \left(\dfrac{1-\xi}{2} \right) \left(\dfrac{1+\xi}{2} \right) 
    P_{p-1}^{(1,1)}(\xi), \quad & 0 < p < P\\[1.5ex]
  \dfrac{1+\xi}{2}, \quad & p=P.
\end{cases}
\label{equ: spectral/hp expansion set}
\end{align}
Here, $\phi_0$ and $\phi_P$ are the linear finite element basis functions and 
\begin{align*}
\phi_1(\xi)=\left(\frac{1-\xi}{2} \right) \left(\frac{1+\xi}{2}\right)
\end{align*}
is the usual quadratic hierarchical expansion mode for quadratic elements.
Furthermore,  $P$ denotes the highest polynomial order of the hierarchical 
expansion and $P_{p}^{(\alpha,\beta)}(\xi)$ denotes the $p^{\mathrm{th}}$-order
Jacobi polynomial.

Two-dimensional functions  $u(\mathbf{x},t)$ can be approximated on two-
dimensional standard quadrilaterals, defined as
$\Omega_{st}=\{-1 \leq\xi_1,\xi_2 \leq 1\}$, using a tensor product of the 
one-dimensional modal expansion basis functions $\phi_p$, that is,
\begin{align}
u(\mathbf{x},t) 
    = \sum_{p=0}^P \sum_{q=0}^P \hat{u}_{pq}(t) \phi_p(\xi_1) \phi_q(\xi_2)
\end{align} 
with the reference coordinates given by
\begin{align}
\xi_1=\left[\chi_1^{e,t}\right]^{-1}(x,y), 
    \quad  \xi_2=\left[\chi_2^{e,t}\right]^{-1}(x,y),
\end{align}  
involving the inverse of the mapping $\pmb{\chi}^e$. Here, the mapping, 
$\pmb{\chi}^e$,  between the local coordinates $(\xi_1,\xi_2)$ and the physical
coordinates $(x,y)$ approximates the geometry with the same order polynomial
space as the solution, that is,
\begin{align}
  \mathbf{x} = \pmb{\chi}^e(\xi_1,\xi_2) = \sum_{p=0}^{P} \sum_{q=0}^{P}
 \mathbf{\hat{x}}_{pq} \phi_{p}(\xi_1) \phi_{q}(\xi_2).
  \label{equ: mappingexpansion}
\end{align}
Details on the construction of this mapping can be found in \cite{Sherwin:2005}.

\subsection{Weak Formulation}
Introducing the function spaces in the current frame with respect to the
reference configuration $\hat{\Omega}_{t_0}$
\begin{subequations}
\begin{align}
  \mathcal{V}(\Omega_t) =& \left\{ \vel: \Omega_t \times I \rightarrow \mathbb{R}^d: \,\, \vel = \hat{\vel} \circ \mathcal{R}_{t}^{-1}, \,\, \hat{\vel} \in   [H^1_0(\hat{\Omega}_{t_0})]^d \right\},  \\
  \mathcal{V}_D(\Omega_t) =& \left\{ \vel: \Omega_t \times I \rightarrow \mathbb{R}^d: \,\, \vel = \hat{\vel} \circ \mathcal{R}_{t}^{-1}, \,\, \hat{\vel} \in   [H^1_D(\hat{\Omega}_{t_0})]^d \right\},  \\
  \mathcal{Q}(\Omega_t) =& \left\{ q: \Omega_t \times I \rightarrow \mathbb{R}^d: \,\, q = \hat{q} \circ \mathcal{R}_{t}^{-1}, \,\, \hat{q} \in  L^2(\hat{\Omega}_{t_0})\right\}, \\
  \mathcal{Q}_0(\Omega_t) =& \left\{ q: \Omega_t \times I \rightarrow \mathbb{R}^d: \,\, q = \hat{q} \circ \mathcal{R}_{t}^{-1}, \,\, \hat{q} \in  L^2_0(\hat{\Omega}_{t_0})\right\}, 
\end{align}
\end{subequations}
the weak formulation of the system of equations~
\eqref{equ: conservation equations ALE-form} leads to the following problem
definition.
\begin{Problem}[Weak formulation of incompressible Navier-Stokes equations]
\label{pro: weak form navier stokes}
For almost every $t \in I$ find $t \rightarrow \left(\vel(t),p(t)\right) \in  \mathcal{V}_D(\Omega_t) \times  \mathcal{Q}_0(\Omega_t)$ such that, for all $\left(\phi_{\vel},\psi\right) \in \mathcal{V}(\Omega_t) \times  \mathcal{Q}(\Omega_t)$
\begin{align}
\Rey \left( \left. \dfrac{\partial \vel}{\partial t} \right|_{\mathbf{Y}} 
+ \left(\vel - \mathbf{w}\right) \cdot \nabla_{\pmb{x}} \vel, \, \phi_{\vel} \right)_{\Omega_t} 
+  \left( 2\rod_{\pmb{x}} , \, \nabla_{\pmb{x}} \phi_{\vel} \right)_{\Omega_t}
\nonumber\\
- \left( p , \, \nabla_{\pmb{x}} \cdot \phi_{\vel}\right)_{\Omega_t} 
- \left\langle  \cauchystress \cdot \mathbf{n}  , \, \phi_{\vel}
\right\rangle_{\Gamma_N(t)} -  \left\langle  \sigma \kappa \cdot \mathbf{n}  ,
\, \phi_{\vel} \right\rangle_{\Gamma_f(t)}&=0,
  \label{equ: ALE momentum}\\
\left(\nabla_{\pmb{x}} \cdot \vel ,\, \psi \right)_{\Omega_t} &= 0 , 
  \label{equ: Weak Form DEVSS-G ALE}
\end{align}
where $\Gamma_N(t)$ is the Neumann boundary and $\Gamma_f(t)$ is the free surface 
boundary.
\end{Problem}

We choose the same trial and test function space for the mesh velocity as for
the fluid velocity, i.e. we choose
\begin{align}
\mathcal{W}(\Omega_t) \equiv \mathcal{V}(\Omega_t)
\end{align}
and we solve Equation~\eqref{equ: mesh velocity strong form} with the boundary 
conditions~\eqref{equ: mesh velocity boundary conditions} using a continuous 
Galerkin method. The weak formulation for the mesh movement can therefore be
expressed as below.
\begin{Problem}[Weak Formulation Mesh Velocity]
\label{pro: weak form mesh velocity}
For almost every $t \in I$ find $t \rightarrow \mathbf{w}(t) \in  \mathcal{V}_D(\Omega_t)$ such that, for all $\phi_{\mathbf{w}} \in \mathcal{V}(\Omega_t)$
\begin{align}
\left(\nabla \mathbf{w}, \nabla \phi_{\mathbf{w}}\right)_{\Omega_t} = 0
\label{equ: elliptic mesh vel}
\end{align}
subject to the boundary conditions~\eqref{equ: mesh velocity boundary conditions}.
\end{Problem}

The position of the new nodes of the mesh can be obtained via 
Equation~\eqref{equ: mesh velocityrepeat}, that is,
\begin{align}
\left.\dfrac{\partial \pmb{x}(\mathbf{Y},t)}{\partial t}\right|_{\mathbf{Y}} 
= \left.\dfrac{\partial \mathcal{R}_{t}(\mathbf{Y}) }{\partial t} 
        \right|_{\mathbf{Y}} = \mathbf{w}(\pmb{x},t). 
\end{align}

\subsection{Discrete ALE formulation}
As mentioned above, we have two referential domains to consider in the ALE
formulation. Firstly, let $\Omega_t^{\delta}$ be the union of all
non-overlapping mesh elements in the Eulerian frame at time $t$ and secondly, let
$\hat{\Omega}_{t_0}^{\delta}$  denote the union of all mesh elements in the
referential frame.
Consider the following discrete trial and test function spaces
\begin{align}
\mathcal{V}_D^{\delta}(\Omega_t^{\delta}) = \left\{ \vel: \Omega_t^{\delta} \times I \rightarrow \mathbb{R}^d: \,\, \vel = \hat{\vel} \circ [\mathcal{R}_t^{\delta}]^{-1}, \,\, \hat{\vel} \in  [H^1_D(\Omega_{t_0}^{\delta})]^d\cap [\mathcal{P}_P^c(\Omega_{t_0}^{\delta})]^d 
\right\}
\end{align}
for the fluid and mesh velocities and
\begin{align}
\mathcal{Q}^{\delta}(\Omega_t^{\delta}) = \left\{ q: \Omega_t^{\delta} \times I \rightarrow \mathbb{R}: \,\, q = \hat{q} \circ [\mathcal{R}_t^{\delta}]^{-1}, \,\, \hat{q} \in  L^2(\Omega_{t_0}^{\delta})\cap [\mathcal{P}_{P-2}(\Omega_{t_0}^{\delta})]^d) 
\right\}, 
\end{align}
for the pressure field.
Alternatively, these spaces can be expressed as (see \cite{Pena2009})
\begin{align}
\mathcal{V}_D^{\delta}(\Omega_t^{\delta}) &= [H^1_D(\Omega_{t}^{\delta})]^d\cap [\mathcal{P}_P^c(\Omega_{t}^{\delta})]^d, \\
\mathcal{Q}^{\delta}(\Omega_t^{\delta}) &=L^2(\Omega_{t}^{\delta})\cap [\mathcal{P}_{P-2}(\Omega_{t}^{\delta})]^d. 
\end{align}
Here, $\mathcal{P}_P^c(\Omega_{t_0}^{\delta})$ denotes the globally continuous space of polynomials of degree $P$ over the reference mesh, that is,
\begin{align}
\mathcal{P}_P^c(\Omega_{t_0}^{\delta}) = \left\{ \left. g^{\delta}: \Omega_{t_0}^{\delta} \rightarrow \mathbb{R} \right| \,\, g^{\delta} \in \mathcal{C}^0(\overline{\Omega_{t_0}}), \,\, \left. g^{\delta} \right|_{\Omega_{t_0}^{e}} \circ [\pmb{\chi}^e(t_0)]^{-1} \in \mathcal{P}_P(\Omega_{st}) \right\}.
\end{align}
$\mathcal{P}_P(\Omega_{t_0}^{\delta})$  denotes the space of piecewise
continuous polynomials of degree $P$ over the reference mesh, that is,
\begin{align}
\mathcal{P}_P(\Omega_{t_0}^{\delta}) = \left\{ \left. g^{\delta}: \Omega_{t_0}^{\delta} \rightarrow \mathbb{R} \right| \,\, g^{\delta} \in L^2(\overline{\Omega_{t_0}}), \,\, \left. g^{\delta} \right|_{\Omega_{t_0}^{e}} \circ [\pmb{\chi}^e(t_0)]^{-1} \in \mathcal{P}_P(\Omega_{st}) \right\}.
\end{align}
$\mathcal{P}_P^c(\Omega_{t}^{\delta})$ denotes the globally continuous polynomial space over the Eulerian mesh  and $\mathcal{P}_P(\Omega_{t}^{\delta}) $ denotes the piecewise continuous polynomial space over the Eulerian mesh. Here, $\left. g^{\delta} \right|_{\Omega_{t_0}^{e}}$ denotes the restriction of $g^{\delta}$ to the spectral element $\Omega_{t_0}^{e}$, $\mathcal{P}_P(\Omega_{st})$ is the space of polynomials of degree $P$ defined on the standard element given by the expansion basis \eqref{equ: spectral/hp expansion set}. Note that, the pressure is discretised with polynomials of order 2 lower than the velocities to satisfy the LBB condition \citep{Brezzi1974}.
The spaces $\mathcal{V}_D^{\delta}(\Omega_t^{\delta})$ and  $\mathcal{Q}^{\delta}(\Omega_t^{\delta})$ include the discrete ALE mapping, which can be expressed as \citep{Nobile2001}
\begin{align}
\left. \mathcal{R}_t^{\delta} \right|_{\Omega_{t_0}^{e}} = \pmb{\chi}^e(t) \circ [\pmb{\chi}^e(t_0)]^{-1} \quad \forall \Omega_{t_0}^{e},
\end{align}
involving the geometrical mappings, $ \pmb{\chi}^e(t)$, at each $t$, from the standard element $\Omega_{st}$ to each element $\Omega_t^e$, that is, 
\begin{align}
\mathbf{x}(\xi_1,\xi_2) =  \pmb{\chi}^e(t;\xi_1,\xi_2) = \sum_{p=0}^{P} \sum_{q=0}^{P} \mathbf{\hat{x}}_{pq}(t) \phi_{p}(\xi_1) \phi_{q}(\xi_2),
\end{align}
where $\mathbf{\hat{x}}_{pq}(t)$ denotes the expansion coefficients at time $t$ and the iso-parametric mapping,  $\pmb{\chi}^e(t_0)$,  from  $\Omega_{st}$ to $\Omega_{t_0}^{e}$,  defined as 
\begin{align}
&\mathbf{Y}(\xi_1,\xi_2) =  \pmb{\chi}^e(t_0;\xi_1,\xi_2) = \sum_{p=0}^{P} \sum_{q=0}^{P} \mathbf{\hat{Y}}_{pq} \phi_{p}(\xi_1) \phi_{q}(\xi_2).
\end{align}


Using these space definitions and an implicit Euler time-integration scheme, the
semi-discrete Navier-Stokes equations are expressed as follows.

\begin{Problem}[Semi-discrete Navier-Stokes ALE formulation]
\label{pro: DEVSSG/DG/AKLE algorithm for vel}
For each $n$, let $t_n=t_0 + n \Delta t$, find
$(\vel^{n+1}_{\delta},p^{n+1}_{\delta})  \in
(\mathcal{V}^{\delta}_D(\Omega_{t_{n+1}}^{\delta}) \times
\mathcal{Q}_0^{\delta}(\Omega_{t_{n+1}}^{\delta}))$ with
$\vel_{\delta}^{0}=\vel_{0,\delta}$ in $\Omega_{t_{0}}^{\delta}$ such that
\begin{eqnarray}
&&\begin{split}
&\Rey \left( \dfrac{\vel^{n+1}_{\delta}-\vel^{n}_{\delta}}{\Delta t}, \, \phi_{\vel} \right)_{\Omega_{t_{n+1}}^{\delta}} 
+ \left( [\left(\vel^*_{\delta} - \mathbf{w}^{n+1}_{\delta}\right) \cdot \nabla_{\pmb{x}}] \vel^{n+1}_{\delta}, \, \phi_{\vel} \right)_{\Omega_{t_{n+1}}^{\delta}} \\
&+  \left( 2\rod^{n+1}_{\pmb{x},\delta} , \, \nabla_{\pmb{x}} \phi_{\vel} \right)_{\Omega_{t_{n+1}}^{\delta}} 
- \left( p^{n+1}_{\delta} , \, \nabla_{\pmb{x}} \cdot \phi_{\vel}\right)_{\Omega_{t_{n+1}}^{\delta}} \\ 
&-  \left\langle  \cauchystress^{n+1}_{\delta} \cdot \mathbf{n}  , \, \phi_{\vel} \right\rangle_{\Gamma_N(t_{n+1})} 
-  \left\langle  \sigma \kappa_S \cdot \mathbf{n}_S  , \, \phi_{\vel} \right\rangle_{\Gamma_f(t_{n+1})}=0, 
 \end{split}
 \label{equ: ALE  momentum discrete}\\
&& \left(\nabla_{\pmb{x}} \cdot \vel^{n+1}_{\delta} ,\, \psi \right)_{\Omega_{t_{n+1}}^{\delta}} = 0 , 
 \label{equ: ALE continuity discrete}
\end{eqnarray}
for all $\left(\phi_{\vel},\psi \right) \in  (\mathcal{V}^{\delta}(\Omega_{t_{n+1}}^{\delta}) \times  \mathcal{Q}^{\delta}(\Omega_{t_{n+1}}^{\delta}))$. 
Here, we linearise the convective term in the momentum equation by setting $\vel^*_{\delta}=\vel^n_{\delta}$, which is an extrapolation of the velocity of the same order as the implicit Euler scheme. Note that, the index $S$ for normals and curvature in the boundary integral over $\Gamma_f(t_{n+1})$ indicates that these quantities are determined from a cubic spline representation of the free surface according to equation \eqref{equ: normal from spline} and \eqref{equ: curv from spline} defined below.
\end{Problem}

\subsection{Matrix formulation}
The discrete ALE formulation involves the following matrices
\begin{eqnarray}
\mathbf{M}^{e}(t)[j][i] &=&   \dfrac{\Rey}{\Delta t}\left( \phi_{u}^{i}, \phi_{u}^{j}, \right)_{{\Omega_{t}}^{e,\delta}}^{\delta}, \\
\mathbf{K}^{e}(t)[j][i] &=& 
\left( \nabla_{\pmb{x}}\phi^i_{\vel} + [\nabla_{\pmb{x}}\phi^i_{\vel}]^T , \, \nabla_{\pmb{x}} \phi^j_{\vel} \right)_{{\Omega_{t}}^{e,\delta}}^{\delta},\nonumber \\
&&- \left\langle  \left(\nabla_{\pmb{x}}\phi^i_{\vel} + [\nabla_{\pmb{x}}\phi^i_{\vel}]^T\right) \cdot \mathbf{n}  , \, \phi^j_{\vel} \right\rangle_{\Gamma_N(t)}, 
 \label{equ: K(t) def}\\
\mathbf{B}^e(t;\vel_{\delta},\mathbf{w}_{\delta})[j][i] &=& \left( [\left(\vel_{\delta} - \mathbf{w}_{\delta}\right) \cdot \nabla_{\pmb{x}}] \phi_{\vel}^{i}, \, \phi^{j}_{\vel} \right)_{{\Omega_{t}}^{e,\delta}}^{\delta}, \\
\mathbf{D}^e(t)[j][i]&=&\left(\nabla_{\pmb{x}} \phi_u^{i}, \psi^j \right)_{{\Omega_{t}}^{e,\delta}}^{\delta},\\ 
\mathbf{b}(t)[j] &=&   \left\langle  \sigma \kappa_S \cdot \mathbf{n}_S  , \, \phi^{j}_{\vel} \right\rangle_{\Gamma_f(t)} ,
\end{eqnarray}
and a modified Helmholtz matrix 
\begin{equation}
\mathbf{H}^e(t)[j][i] := \mathbf{M}^{e}(t)[j][i] +  \mathbf{K}^{e}(t)[j][i] + \mathbf{B}^e(t;\vel_{\delta},\mathbf{w}_{\delta})[j][i].
\end{equation}
The equation system \eqref{equ: ALE  momentum discrete}-\eqref{equ: ALE continuity discrete}  can then be written for each element in algebraic form as 
\begin{align}
\mathbf{H}_{g}(t_{n+1}) \hat{\vel}^{n+1}_g - \mathbf{D}_{g}(t_{n+1})^{T} \hat{\mathbf{p}}^{n+1}_g &= 
\mathbf{M}(t_{n+1}) \hat{\vel}^{n} + \mathbf{b}(t_{n+1}), \nonumber\\
\mathbf{D}_{g}(t_{n+1})\textit{} \hat{\vel}_g^{n+1} &= 0, 
\label{equ: coupled system upG}
\end{align}
where $\hat{\vel}_g$ and $\hat{\mathbf{p}}_g$ are the vectors of unknown global coefficients, $\mathbf{H}_g$, $\mathbf{D}_g=(\mathbf{D}_{x_1}, \mathbf{D}_{x_2})$ are the global matrices assembled from the elemental matrix contributions. 
The resulting system of equations is then solved using a multi-level static
condensation technique introduced by \cite{Ainsworth1999}, \cite{Sherwin2000}
and \cite{Sherwin:2005} for the Stokes equations in fixed domains.
 
\subsection{Discretisation of Mesh Movement}
Even though solving Problem~\ref{pro: weak form mesh velocity} yields
continuous mesh movement, the free surface boundary might not be sufficiently
smooth. The free surface boundary undergoes the largest deformation and its
movement involves the evaluation of outward normals, $\mathbf{n}$, in
Equation~\eqref{equ: mesh velocity boundary conditions}, across multiple
elements. Note that, a standard Galerkin method with a $C^0$-continuity across
elements is not sufficient to determine a well-defined normal at element edges.
To alleviate this problem, we represent the free surface using a cubic spline,
$S(x,t) \in \mathcal{C}^2(\Gamma_f)$ to ensure sufficient smoothness of free
surface boundary edges of the mesh. The cubic spline can then be used to
determine the unit outward normals $\mathbf{n}$ and the curvature $\kappa$ of
the free surface using
\begin{align}
\mathbf{n}_S(t) = \dfrac{1}{\sqrt{S'(x,t)^2+1}} 
\left(
\begin{array}{c}
-S'(x,t) \\ 
1
\end{array} \right), 
\label{equ: normal from spline}\\
\kappa_S(t) = \dfrac{|S''(x,t)|}{(1+S'(x,t)^2)^{3/2}}.
\label{equ: curv from spline}
\end{align}
These expressions are then used to evaluate the free surface boundary condition for the mesh velocity given by Equation~\eqref{equ: mesh velocity boundary conditions} 
and the free surface boundary integral in the momentum equation
\begin{align}
\int\limits_{\Gamma_f} \sigma \kappa_S \mathbf{n}_S \, \phi_{\vel} \,
\mathrm{d}\Gamma.
\end{align}
For given $\vel^n$, we perform the mesh movement in the following way. First, we determine the cubic spline through all the quadrature points along the free surface.
Let $(x_i,y_i)$, $1\leq i \leq N,$ be the physical coordinates of the $N$ quadrature points along the free surface. Then, we construct a cubic spline $S(x,t)=S_i(x,t)$ for each $x_i \leq x \leq x_{i+1}$ through
\begin{align}
S_{i}(x,t) = a_i (x-x_i)^3 + b_i (x-x_i)^2 + c_i (x-x_i) + d_i
\end{align}
where we enforce \emph{continuity} 
\begin{align}
S_{i-1}(x_i,t) = S_i(x_i,t), \nonumber \\
S_{i}(x_{i+1},t) = S_{i+1}(x_{i+1},t)
\end{align}
and \emph{smoothness} 
\begin{align}
S_{i-1}'(x_i,t) = S_i'(x_i,t), \nonumber \\
S_{i-1}''(x_i,t) = S_i''(x_i,t), \nonumber \\
S_{i}'(x_{i+1},t) = S_{i+1}'(x_{i+1},t), \nonumber \\
S_{i}''(x_{i+1},t) = S_{i+1}''(x_{i+1},t).
\end{align}
We employ the \emph{not-a-knot} boundary condition on the spline, that is, 
\begin{align}
S'''_1(x_2)=S'''_{2}(x_2),\\ 
S'''_{N-1}(x_{N-1})=S'''_{N-2}(x_{N-1}).
\end{align} 
We then solve the elliptic problem~\eqref{equ: elliptic mesh vel} using the
continuous Galerkin method, determining
\begin{align}
\mathbf{L}_g \hat{\tilde{\mathbf{w}}}_g = 0,
\label{equ: intermediate mesh vel}
\end{align}
where $\mathbf{L}_g$ is the global Laplace matrix given by  
\begin{align}
\mathbf{L}^{e}(t)[j][i] =	 
\left( \nabla_{\pmb{x}}\phi^i_{\mathbf{w}} , \, \nabla_{\pmb{x}}
\phi^j_{\mathbf{w}} \right)_{{\Omega_{t}}^{e,\delta}}^{\delta},
\end{align}
subject to the boundary conditions~\eqref{equ: mesh velocity boundary conditions}, 
which include the normal determined by the cubic spline according 
to~\eqref{equ: normal from spline}.

The mesh velocity resulting from the solution of Equation~\eqref{equ:
intermediate mesh vel}, denoted by $\tilde{\mathbf{w}}$, is then used to update
the coordinates of the mesh nodes using a third order Adams-Bashforth-Scheme for
Equation~\eqref{equ: mesh velocityrepeat}, that is,
\begin{align}
\mathbf{X}^{n+1} = \mathbf{X}^n + \frac{\Delta t}{12} (23 \tilde{\mathbf{w}} - 16 \mathbf{w}^n + 5 \mathbf{w}^{n-1}).
\end{align}
This equation is solved pointwise in the strong form for each quadrature point. 
However, in practice, we do not move all the mesh nodes of every element. We 
only move all the quadrature points along the free surface boundary introducing
curved edges along the free surface boundary. In the interior of the domain, we
just move the corner vertices of every element keeping the interior edges of
the domain straight.

Using the new coordinates of all mesh nodes, we compute the mesh velocity at 
the new time level pointwise as 
\begin{align}
\mathbf{w}^{n+1} = \dfrac{\mathbf{X}^{n+1}-\mathbf{X}^{n}}{\Delta t}.
\end{align}
\newpage
\subsection{Algorithm Summary}
In Summary, the solution procedure is outlined in Algorithm~\ref{alg: DEVSS-G ALE}
 \begin{center}
 \begin{pseudocode}[shadowbox]{ALE scheme.}{\mathbf{u}^n,p^n}
 \label{alg: DEVSS-G ALE}
 t = t_0\\
 \WHILE t\leq t_{fin} \DO
 \BEGIN
 \PROCEDURE{MoveMesh}{\mathbf{u}^n,p^n,\es^n}
 \text{Construct Cubic Spline through Free Surface Boundary.}\\
 \text{Set BC for Mesh Velocity (see \eqref{equ: mesh velocity boundary conditions}).}\\
 \text{Solve Elliptic Problem for Mesh Velocity \eqref{equ: mesh veldiscrete}.} \\
 \OUTPUT{\mathbf{w}^{n+1}}\\
 \text{Compute New Mesh Coordinates $\mathbf{X}^{n+1}$.}\\
 \text{Construct New Parametric Mappings $\pmb{\chi}^{e}(t_{n+1})$.}\\
 \OUTPUT{\Omega_{t_{n+1}}}
 \ENDPROCEDURE
 \text{Set Boundary Conditions for $\vel$ and $p$.}\\
 \PROCEDURE{SolveCoupledSystem}{\mathbf{u}^n,p^n,\mathbf{w}^{n+1} }
\text{Solve Coupled System of Velocity, Pressure} \\
 \OUTPUT{\mathbf{u}^{n+1},p^{n+1}}
 \ENDPROCEDURE
 t_{n+1} \GETS t_{n} + \Delta t \\
 n+1 \GETS n
 \END\
 \end{pseudocode} 
 \end{center}

\section{Numerical Results}
\label{sec: numerical results}
\subsection{Mesh Configuration}
\begin{figure}[t] \centering
\subfloat{\includegraphics[width=\linewidth]{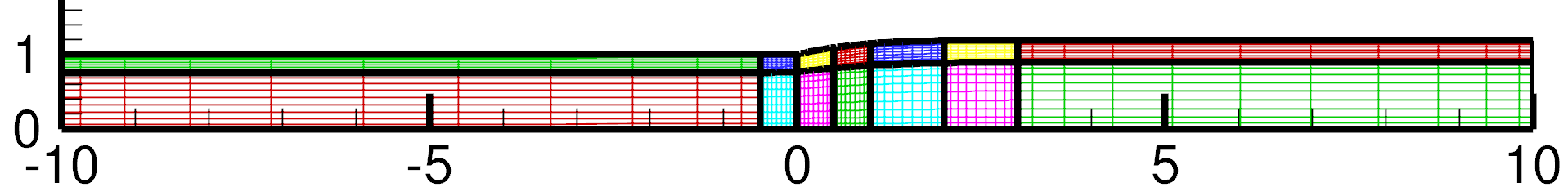}}\\
\subfloat{\includegraphics[width=\linewidth]{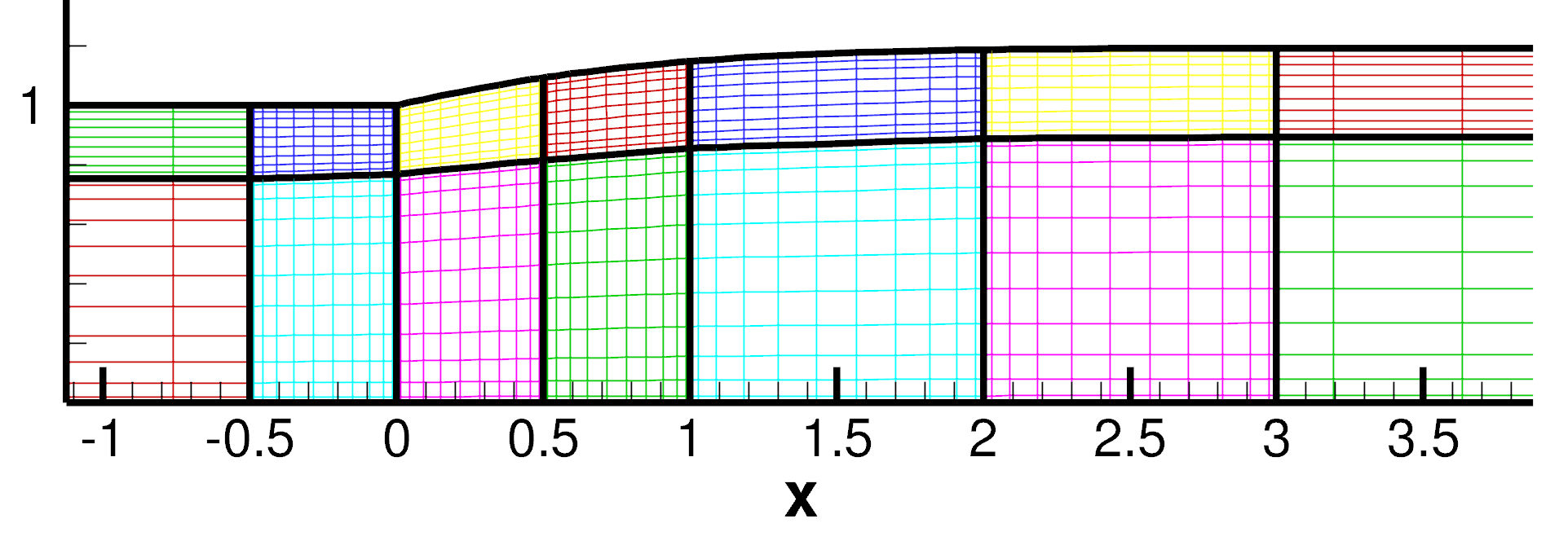}}
\caption{Mesh configuration used for the extrudate swell computation.}
\label{fig:DieSwellGeoNewtonian}
\end{figure}
We use a mesh consisting of $N_{el}=14$ elements as shown in
Figure~\ref{fig:DieSwellGeoNewtonian}  and refine the mesh by increasing the
polynomial order $P$.
We consider a die of length $L_1=10$ and an exit region of length $L_2=10$. The entry length is sufficient to guarantee a fully developed flow far
upstream from the exit of the die. The exit length is chosen sufficiently long to allow the free surface to reach a
constant downstream height for a large range of Reynolds numbers.  For high Reynolds numbers, the free jet length might be insufficient to guarantee a fully developed plug flow profile at outflow. However, the use of the open outflow boundary condition enables us to predict the correct swelling ratios truncated at the outflow boundary location (see  \cite{Mitsoulis2012a}). Throughout this section, we choose a time step of $5\times 10^{-3}$.

\subsection{Numerical results for $Re=0$}

   \begin{table}[t]
    \centering
    \ra{1.2}
    \footnotesize
    \caption{Newtonian swelling ratios for $\Rey=0$}
    \rowcolors{2}{gray!25}{white}
    \begin{tabular}{llll}
    \mytoprule 
     & Method & DOF &  $\chi_R$  \\  
      \mymidrule  
\cite{Crochet1982} & FEM & 562 & 1.200    \\
 &  & 1178 & 1.196    \\
\cite{Reddy1978} & FEM & 254 & 1.199 \\
\cite{Mitsoulis2012a} & FEM & 11270 & 1.191 \\
 &  & 30866 & 1.186 \\
 \cite{Georgiou1999} & FEM (SFEM) & 7528 & 1.1919 (1.1863) \\
  & FEM (SFEM)& 12642 & 1.1888  (1.1863) \\
   \mybottomrule
    \end{tabular}
   \label{tab: Newtonian Swelling}
    \end{table}

Inertialess Newtonian extrudate swell has been investigated in a number of
publications. Table~\ref{tab: Newtonian Swelling} summarises some of the
swelling ratios obtained by a range of authors for plane Newtonian die swell.
\cite{Tanner:2002} used the results in the literature to estimate an
extrapolated value for planar die swell of $\chi_R=1.190\pm 0.002$.  In
general, an increase in the number of degrees of freedom yields less swelling.
  
   \begin{table}[t]
   \centering
   \ra{1.2}
   \footnotesize
   \caption{Comparison of swell ratios and exit pressure corrections for increasing number of degrees for freedom (DOF) between our algorithm and  \cite{Taliadorou2007}.}
   \rowcolors{5}{white}{gray!25}
   \begin{tabular}{l*{7}{l}}
   \mytoprule
   \multicolumn{4}{c}{ Spectral/hp method} & \multicolumn{3}{c}{\cite{Taliadorou2007} FEM} \\
   \cmidrule(r){1-4} \cmidrule(r){5-7}
   P &	DOF & $h_f$ & $n_{ex}$ & DOF &  $h_f$ & $n_{ex}$		\\
   \mymidrule
   8 & 2624	&1.1928 & 0.1507 & & & \\
   10 & 4116	&1.1912 &  0.1503 & 37208 &  1.1953 &	0.1514  \\	
   12 & 5944	&1.1901 & 0.1497 & 43320 &	 1.1908 & 0.1491 \\
   14 & 8108	&1.1900 &  0.1491 & 	49864 &	 1.1893 & 0.1482  \\
   16 & 10608  &1.1891 &  0.1485 & 60490 &  1.1878 &	0.1473	\\
   \mybottomrule
   \end{tabular} 
   \label{tab: pressure exit correction comparison}
   \end{table} 

\begin{figure}[t]
\centering

\subfloat[Velocity component $u$ along the free
surface.]{\includegraphics[width=.5\linewidth]{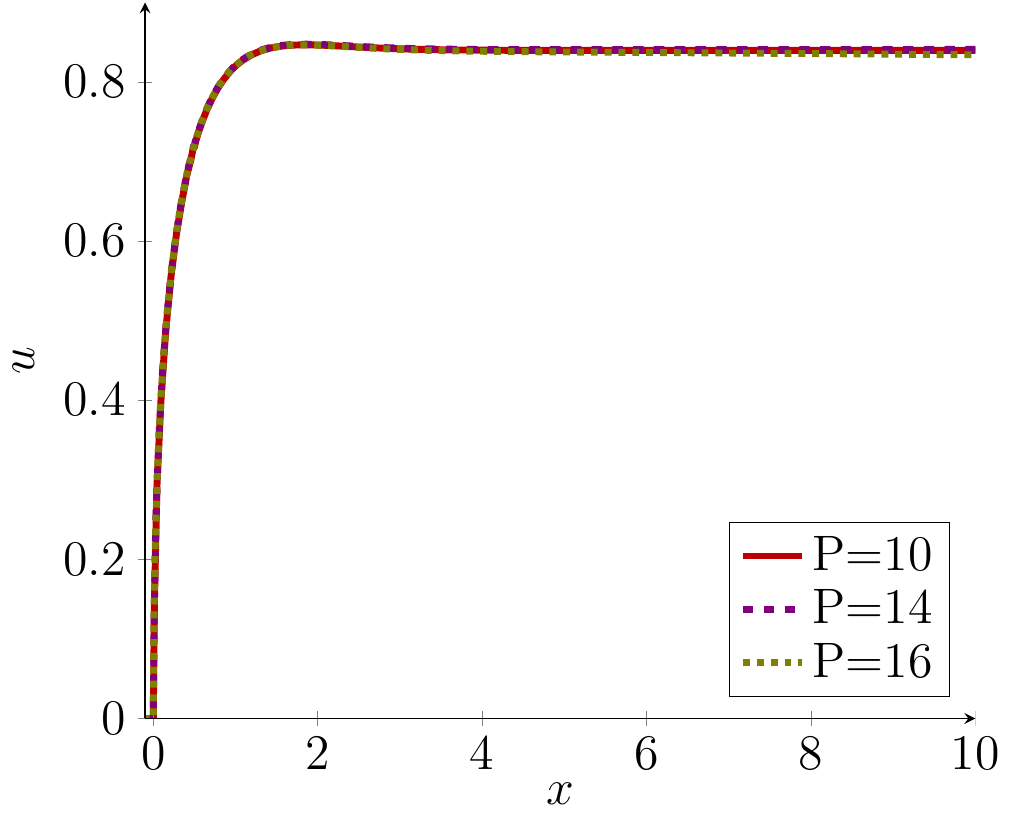}\label{subfig:
swell newton refinement u}} \subfloat[Velocity component $v$ along the free
surface.]{\includegraphics[width=.5\linewidth]{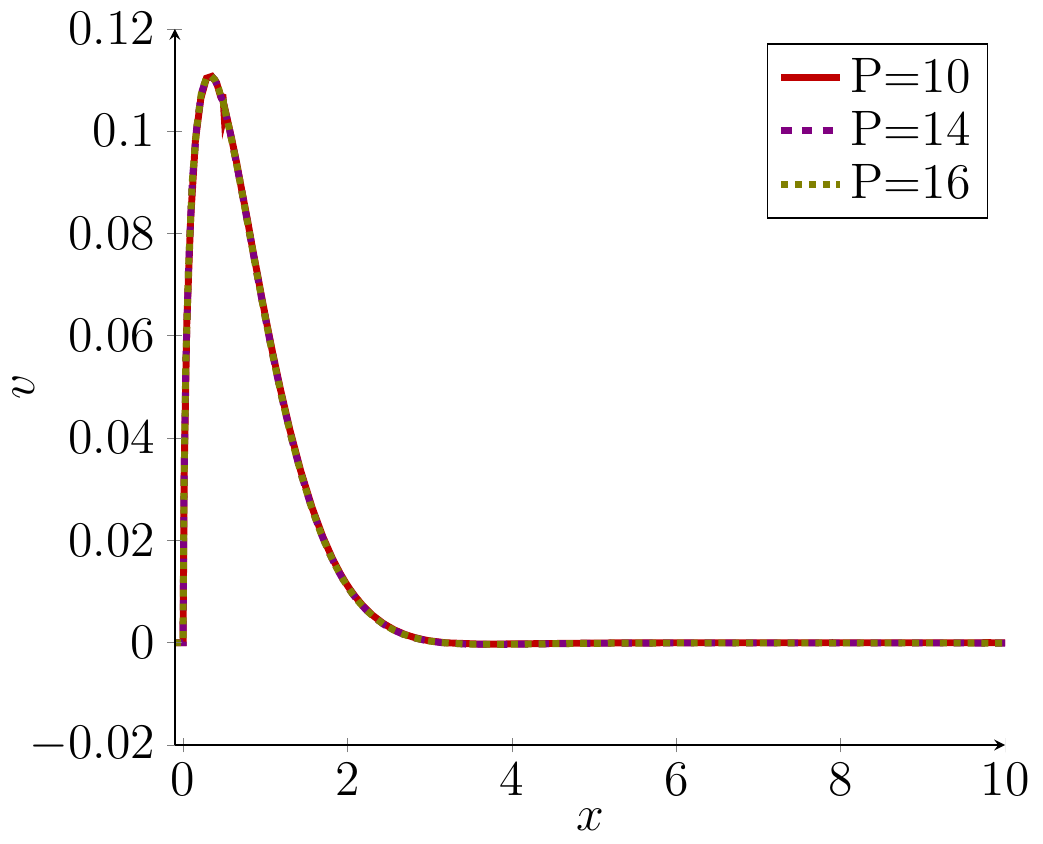}\label{subfig:
swell newton refinement v}}\\
\subfloat[Pressure $p$ along the die wall and the free
surface.]{\includegraphics[width=.65\linewidth]{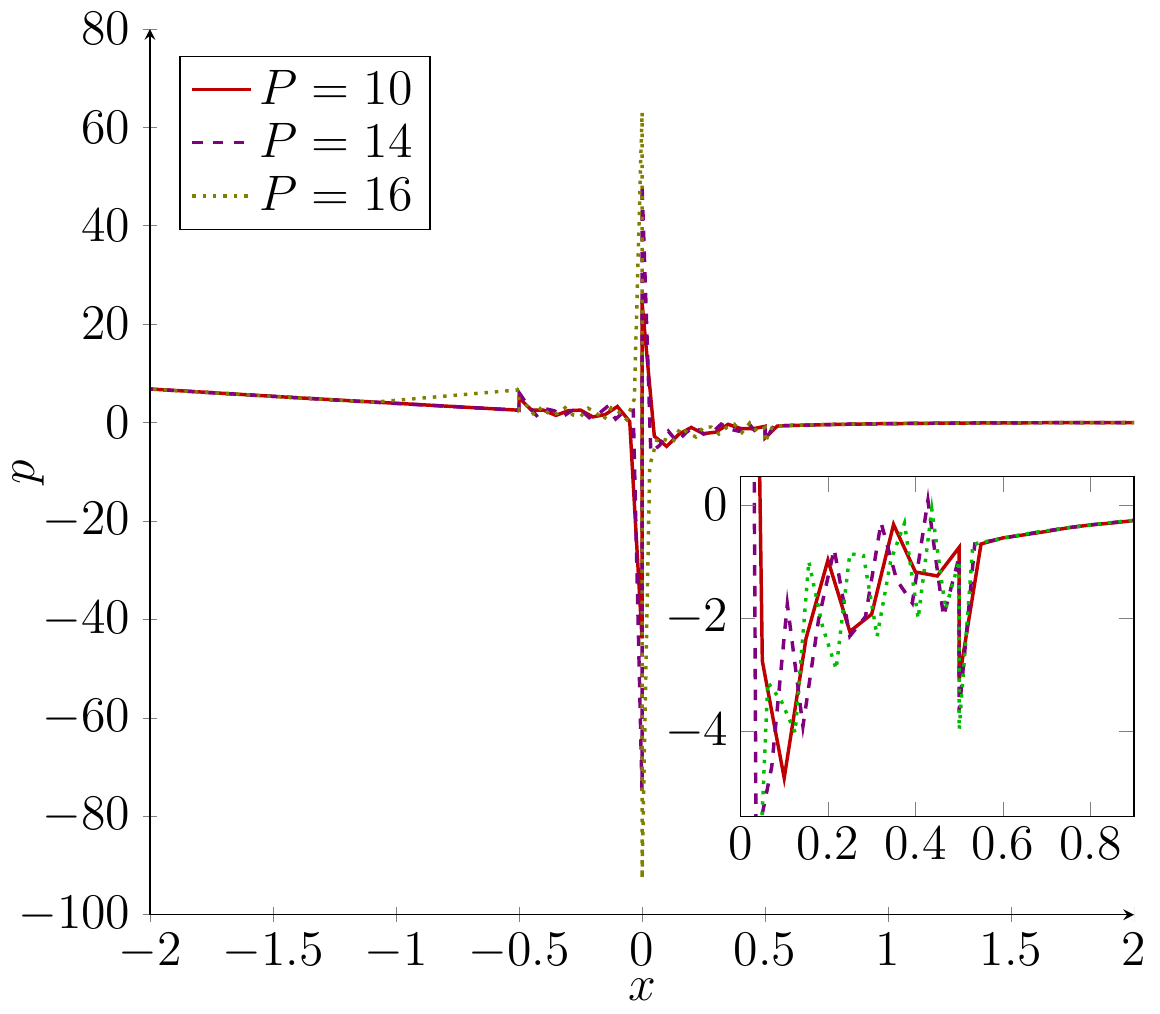}\label{subfig:
swell newton refinement p}} \subfloat[Pressure minima and
maxima.]{\includegraphics[width=.35\linewidth]{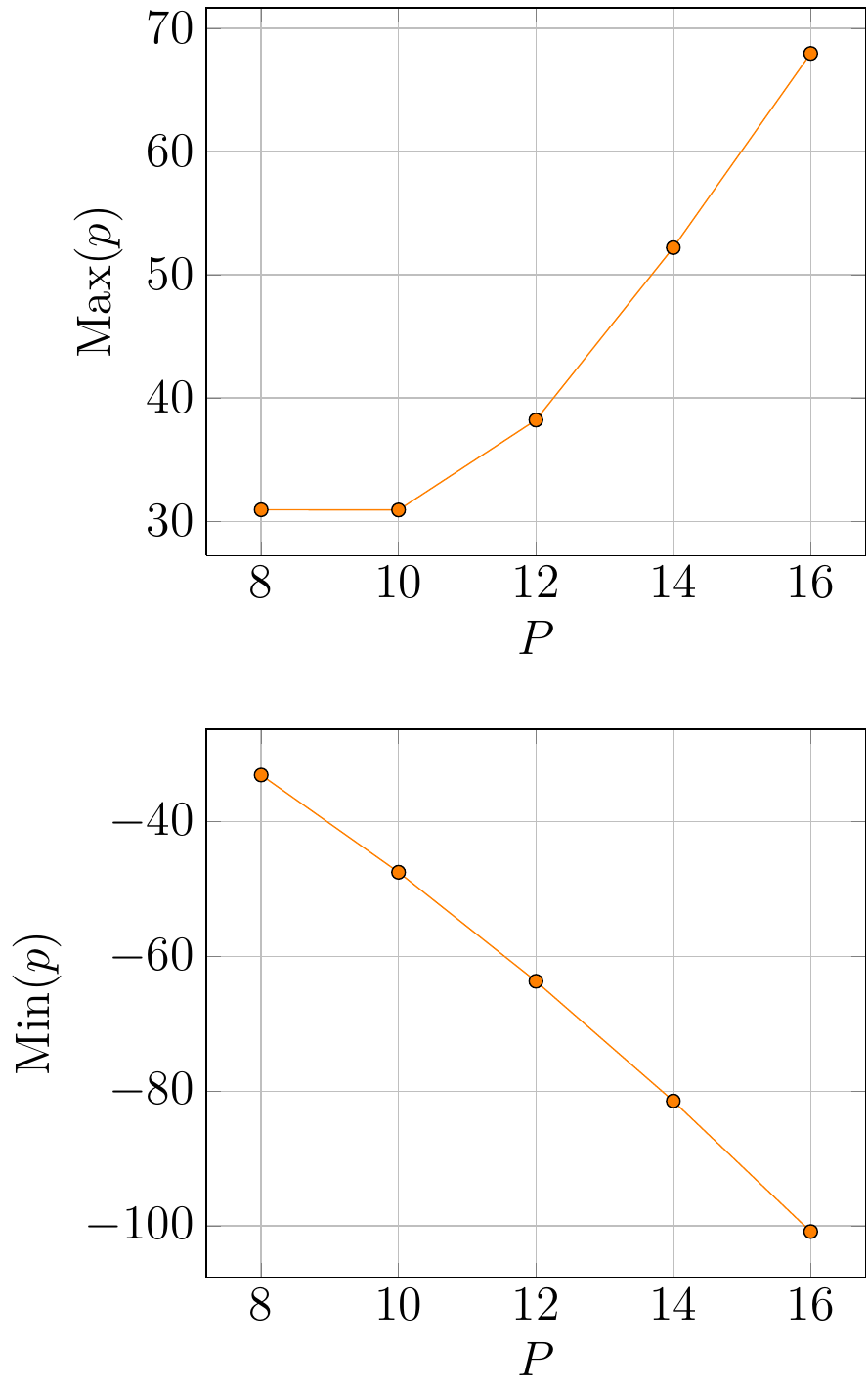}\label{subfig:
p minmax refinement}}\\
\caption{Influence of $P$-mesh refinement on \protect\subref{subfig: swell newton refinement u} the velocity components $u$, \protect\subref{subfig: swell newton refinement v} $v$ and \protect\subref{subfig: swell newton refinement p} pressure $p$ along the free surface and the increase of maximum and minimum values of the pressure at the singularity with increasing polynomial order \protect\subref{subfig: p minmax refinement}.}
\label{fig: swell newton refinement uvp}
\end{figure}

Table~\ref{tab: pressure exit correction comparison} lists a comparison of the
pressure exit correction for $\Rey=0$ of our scheme and the swell ratio for
increasing mesh refinement with the results obtained by \cite{Taliadorou2007}. We
obtain close agreement for a much smaller number of degrees of freedom, which
demonstrates that $p$-refinement is effective for the Newtonian extrudate swell
even though the result is polluted by Gibbs oscillations in the pressure around
the singularity (Figure~\ref{subfig: swell newton refinement p}). The Gibbs
oscillations in the pressure  stay confined to the elements adjacent to the
singularity.    Increasing the Reynolds number leads to a dampening in the
oscillations in the elements adjacent to the singularity and the extreme values
of the pressure at the singularity decrease significantly (Figure~\ref{subref:
fs newt swell p Re}). As shown in Figure~\ref{subfig:
swell newton refinement p} increasing the polynomial order yields an increase in the number of
oscillations. However, the amplitude of each oscillation is reduced with
increasing polynomial order $P$. Increasing the polynomial order also has the
effect of exponentially increasing the maximum value of the pressure and sharply
increasing the minimum value of the pressure at the singularity which reflects
an improved approximation of the infinite pressure value at the singularity
(Figure~\ref{subfig: p minmax refinement}). While the infinite pressure values
at the singularity hamper the rate of convergence of the numerical pressure
solution, the values of the velocity components along the free surface are
converged for $P\geq 10$ (see Figure~\ref{fig: swell newton refinement
uvp}\protect\subref{subfig: swell newton refinement u}, \protect\subref{subfig:
swell newton refinement v}).

\subsection{Impact of inertia}

\begin{figure}[t]
\centering
\includegraphics[width=0.7\linewidth]{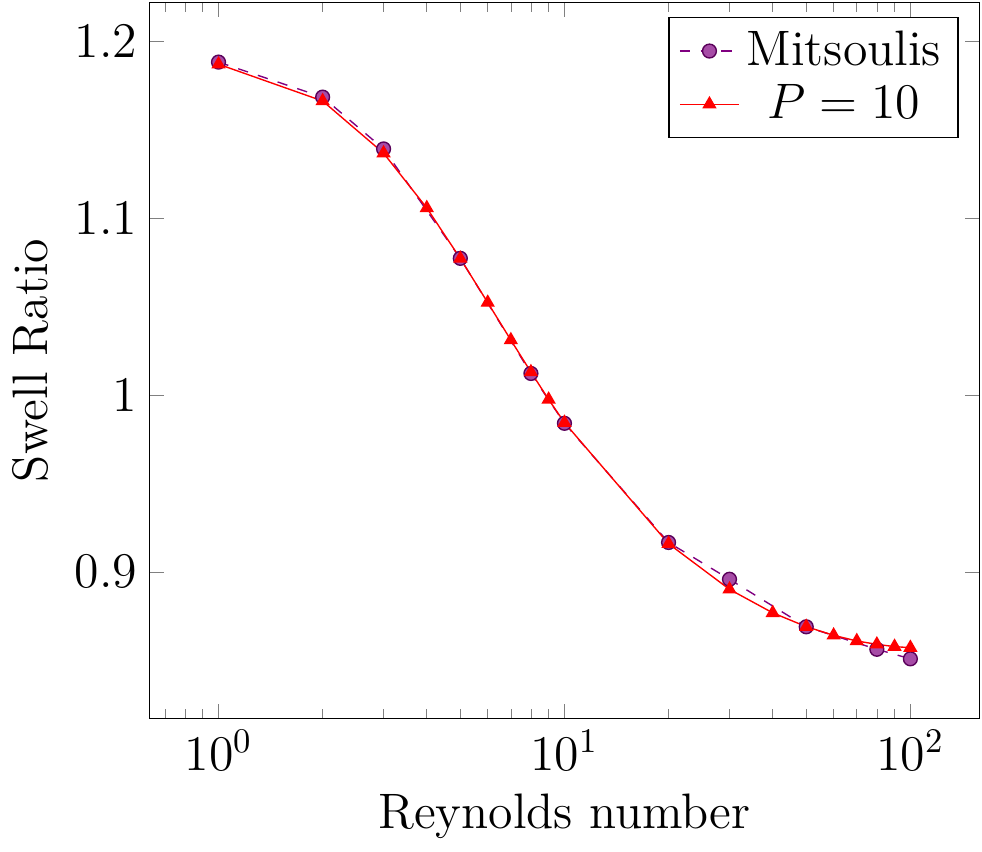}
\caption{Comparison of swell ratios for Newtonian fluid from the current
study ($P=10$) with Mitsoulis \emph{et al}\cite{Mitsoulis2012a}.}
\label{fig:NewtSwellMitsoulis}
\end{figure}

\begin{table}[t]
   \centering
   \ra{1.2}
   \footnotesize
   \caption{Comparison of Newtonian die swell ratio for increasing Reynolds number with \cite{Mitsoulis2012a}.}
   \rowcolors{3}{gray!25}{white}
   \begin{tabular}{lp{2cm}llp{2cm}l}
   \mytoprule 
   $\Rey$  & \cite{Mitsoulis2012a} &$P=10$ &  $\Rey$  & \cite{Mitsoulis2012a} &$P=10$ \\  
     \mymidrule  
   0 & 1.1915 & 1.1912 & 10 & 0.9842 & 0.9846 \\ 
   1 & 1.1885 & 1.1873 & 20 & 0.9168 & 0.9161 \\ 
   2 & 1.1687 & 1.1665 & 30 & 0.8960 & 0.8903 \\ 
   3 & 1.1394 & 1.1370 & 40 & \multicolumn{1}{l}{} & 0.877 \\ 
   4 & \multicolumn{1}{l}{ } & 1.1060 & 50 & 0.8691 & 0.8692 \\ 
   5 & 1.0775 & 1.0774 & 60 & \multicolumn{1}{l}{} & 0.8643 \\ 
   6 & \multicolumn{1}{l}{} & 1.0525 & 70 & \multicolumn{1}{l}{} & 0.8611 \\ 
   7 & \multicolumn{1}{l}{} & 1.0313 & 80 & 0.8564 & 0.8592 \\ 
   8 & 1.0124 & 1.0132 & 90 & \multicolumn{1}{l}{} & 0.8579 \\ 
   9 & \multicolumn{1}{l}{} & 0.9977 & 100 & 0.85103 & 0.8573 \\ 
  \mybottomrule
   \end{tabular}
  \label{tab: Swelling Newton Mitsoulis}
\end{table}

Inertia causes a decrease of the swelling and the liquid jet eventually
contracts for sufficiently high Reynolds numbers. 
We performed computations for Reynolds numbers ranging from 0 to 100. We start by 
computing the extrudate swell for $\Rey=0$ and initialise this
computation with the solution of the corresponding stick-slip problem. After
having obtained the extrudate swell for $\Rey=0$, we increase the
Reynolds number in steps of 1 from 1 to 10 and in steps of 10 from 10 to 100,
each time using the result of the converged extrudate swell of the previous
lower Reynolds number as the initial condition. As the convergence criterion, we
choose a change of the maximum absolute value of all variables including the
mesh velocity of less than $10^{-6}$. Figure~\ref{fig:NewtSwellMitsoulis} and
Table~\ref{tab: Swelling Newton Mitsoulis} shows the comparison of the swelling
ratios obtained with our algorithm with the results of \cite{Mitsoulis2012a},
which are in excellent agreement.

\begin{figure}[t]
\centering
\subfloat{\includegraphics[width=.5\linewidth]{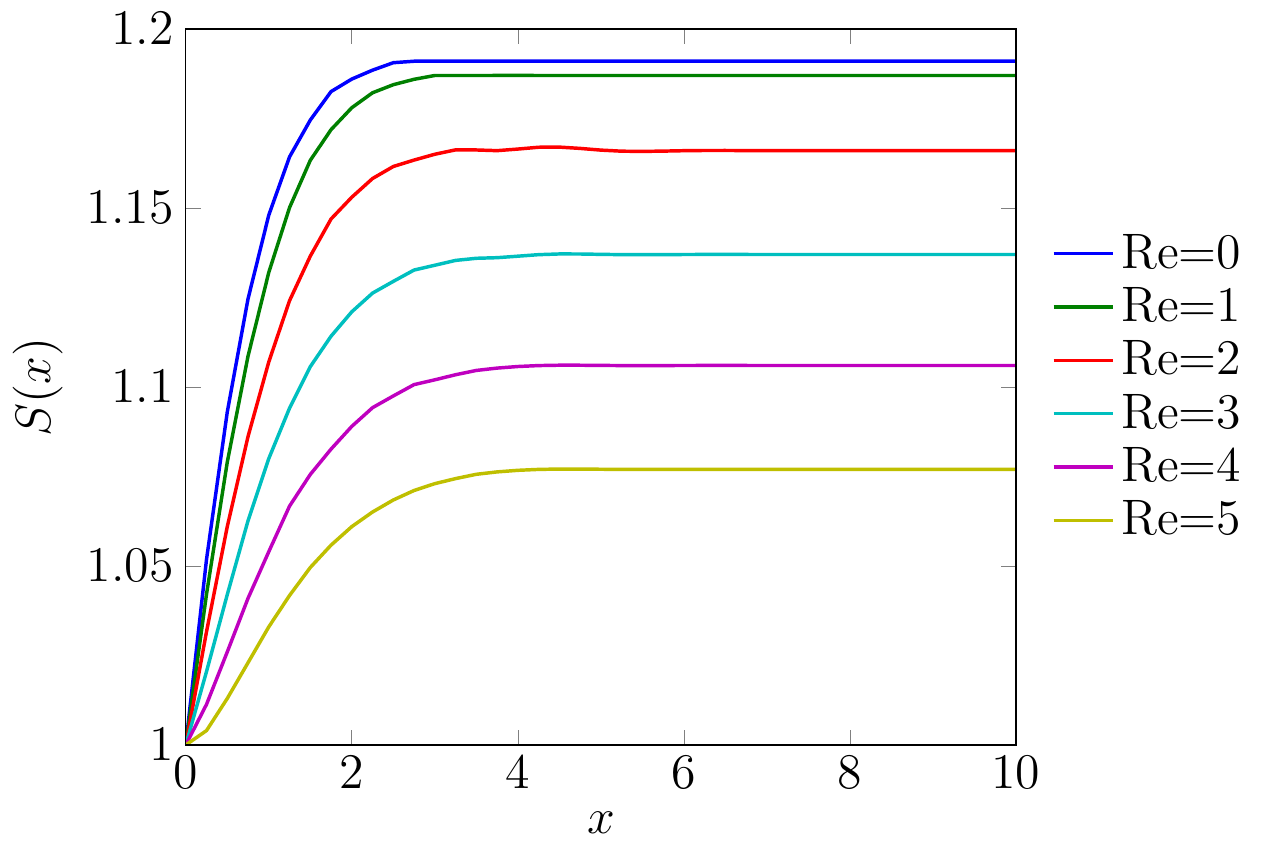}}
\subfloat{\includegraphics[width=.5\linewidth]{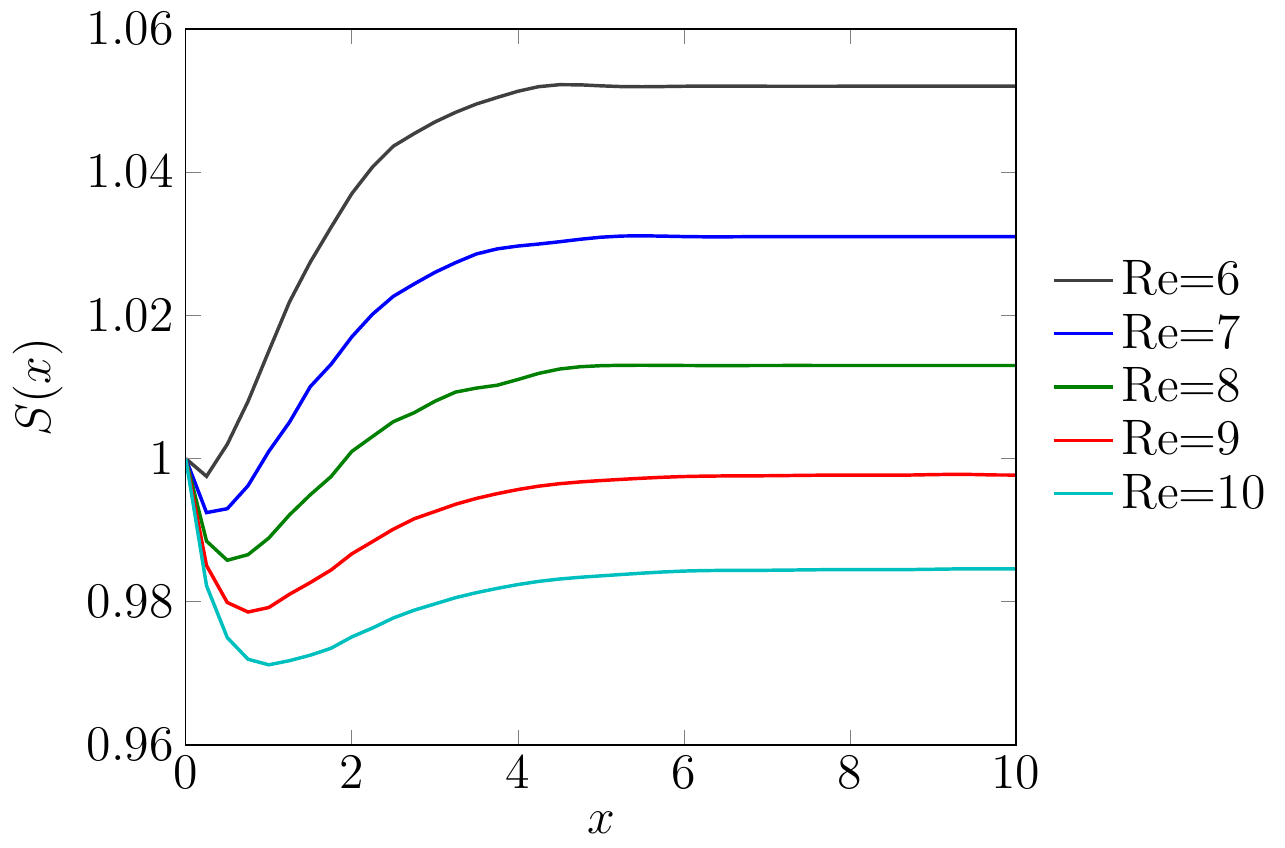}}\\
\subfloat{\includegraphics[width=.5\linewidth]{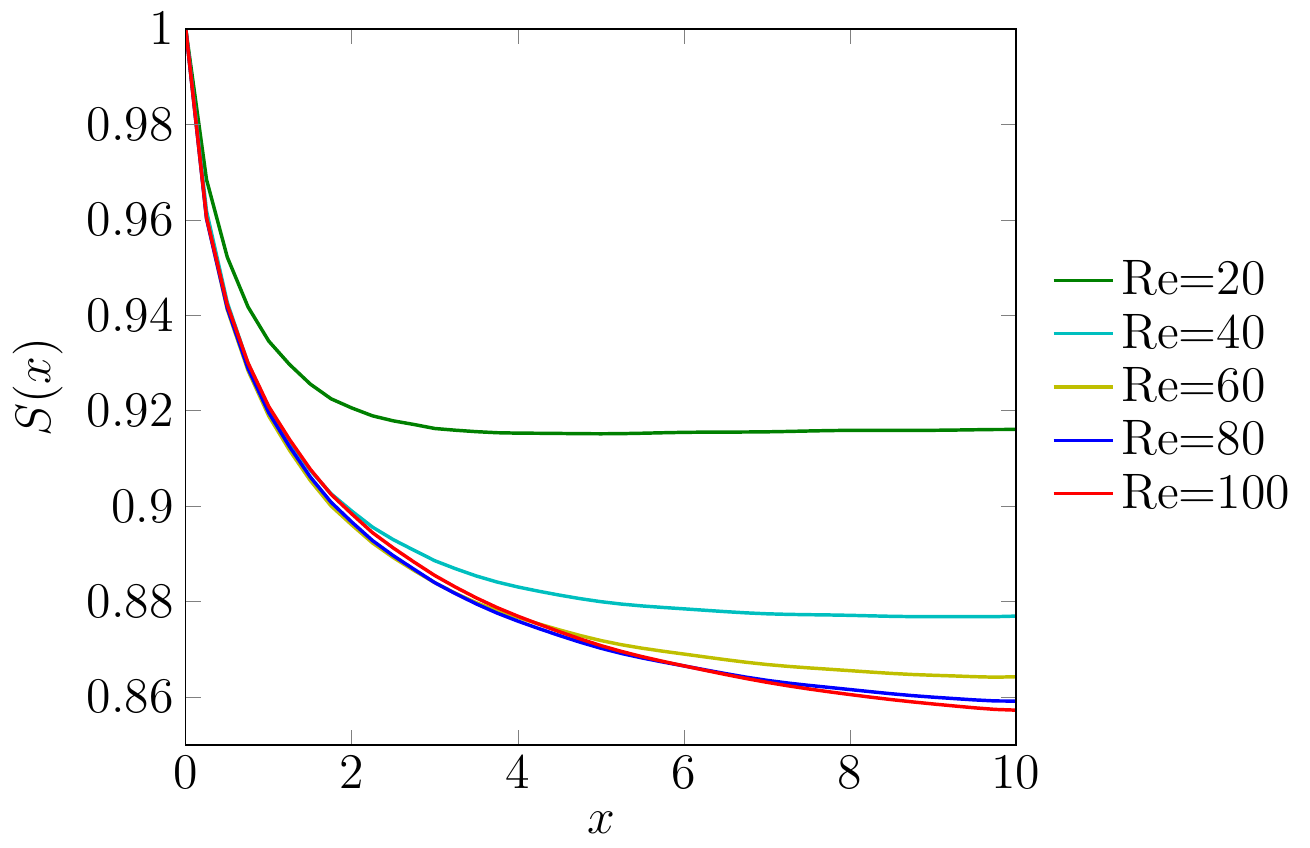}}
\caption{Free surface spline profiles for Newtonian extrudate swell for $P=10$ for a range of Reynolds numbers.}
\label{fig:NewtSwellRe10}
\end{figure}

Figure~\ref{fig:NewtSwellRe10} displays the corresponding free surface spline
profiles.
We observe that the swelling ratio decreases at an accelerating pace with
increasing Reynolds number until $\Rey=6$. For $\Rey=6$, we
see the onset of a delayed die swell in which the fluid surface first goes
through a minimum before it swells again. The delay in the swelling of the jet
increases with increasing Reynolds number from $\Rey=6$ to $\Rey=10$. For
$\Rey=9$ and $\Rey=10$, the fluid contracts ($\chi_R<1$) but still experiences
some swelling after going through a minimum near the die exit. For $\Rey=20$ to
$\Rey=100$ the fluid does not experience any delayed swelling and contracts. For
$10<\Rey<40$ the fluid contracts very fast with increasing Reynolds number. This
trend in the contraction rate with increasing Reynolds number then slows down
and approaches a limit for $40<\Rey<100$. The limit for infinite Reynolds number
was estimated by \cite{Tillett1968} who performed a boundary layer analysis for
a free Newtonian jet and predicted a limiting value of $\chi_R=0.8333$ for
infinite Reynolds number.


\begin{figure}[t]
\centering
\begin{tabular}{ p{0.1\linewidth} p{0.8\linewidth}}
(a) & \raisebox{-0.9\height}{\includegraphics[width=\linewidth]
                    {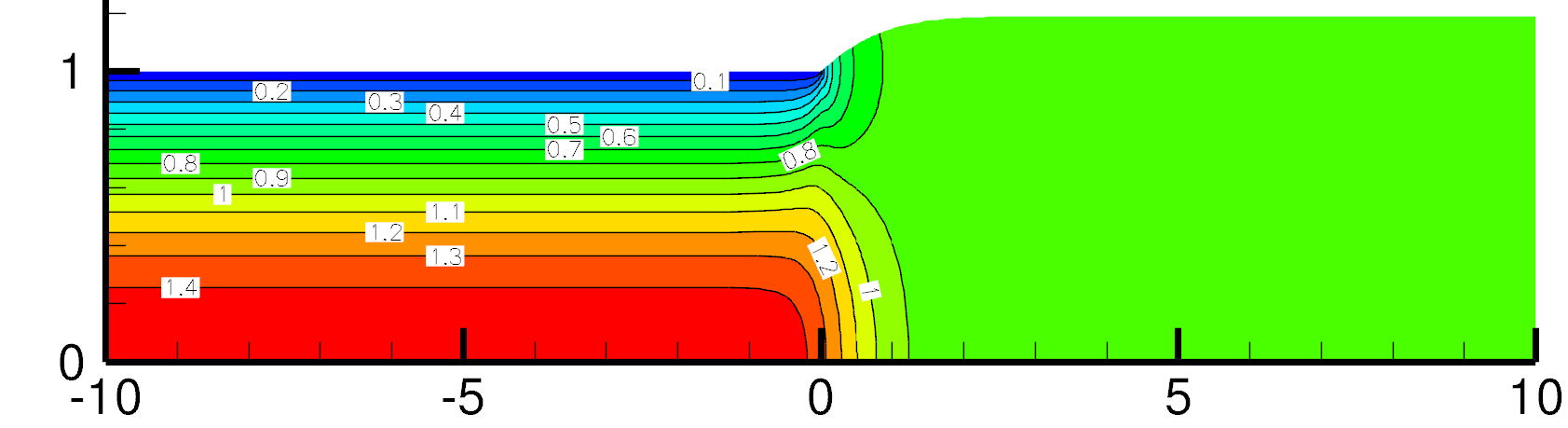}}\\
(b) & \raisebox{-0.9\height}{\includegraphics[width=\linewidth]
                    {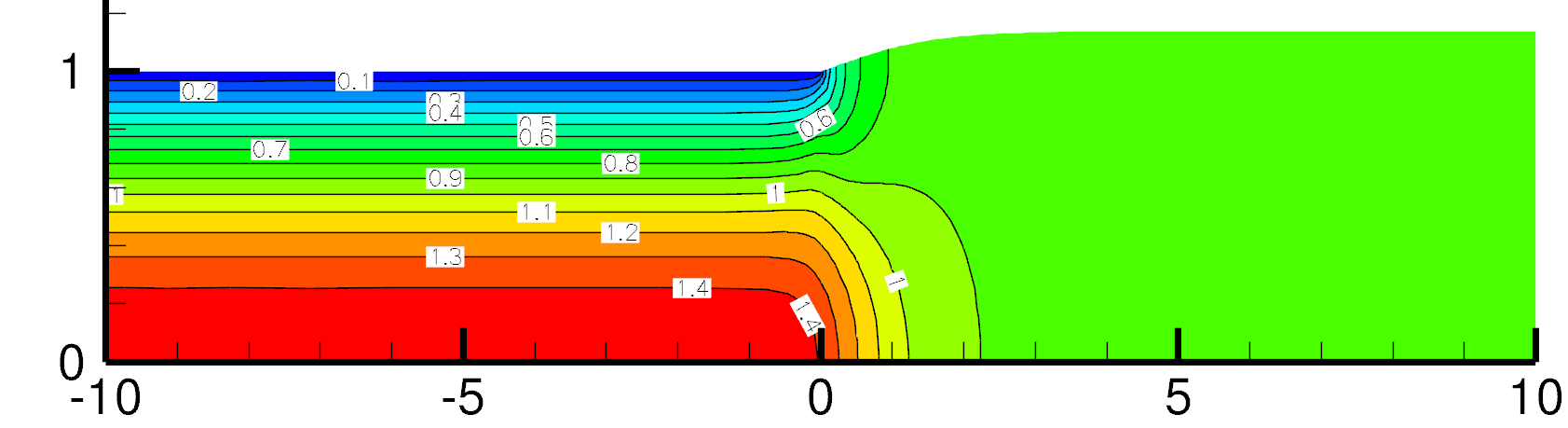}}\\
(c) & \raisebox{-0.9\height}{\includegraphics[width=\linewidth]
                    {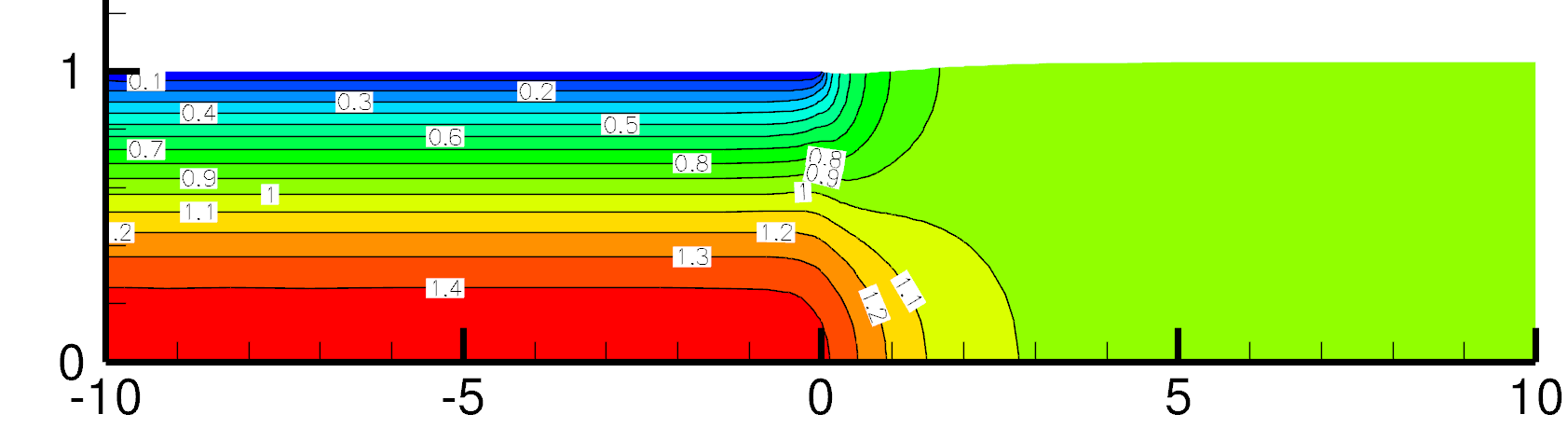}}\\
(d) & \raisebox{-0.9\height}{\includegraphics[width=\linewidth]
                    {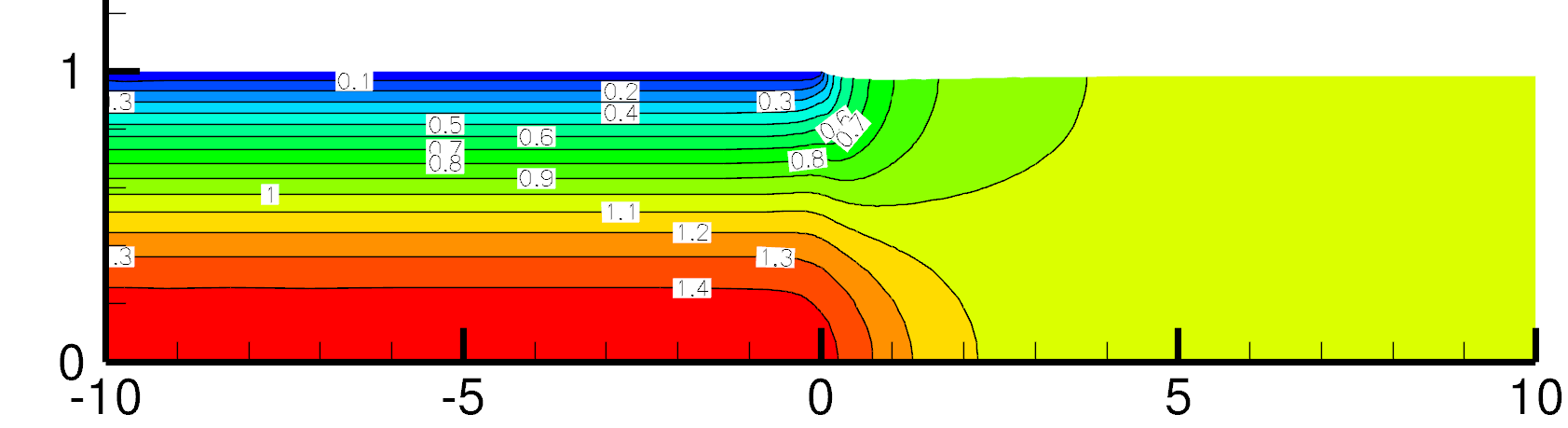}}\\
(e) & \raisebox{-0.9\height}{\includegraphics[width=\linewidth]
                    {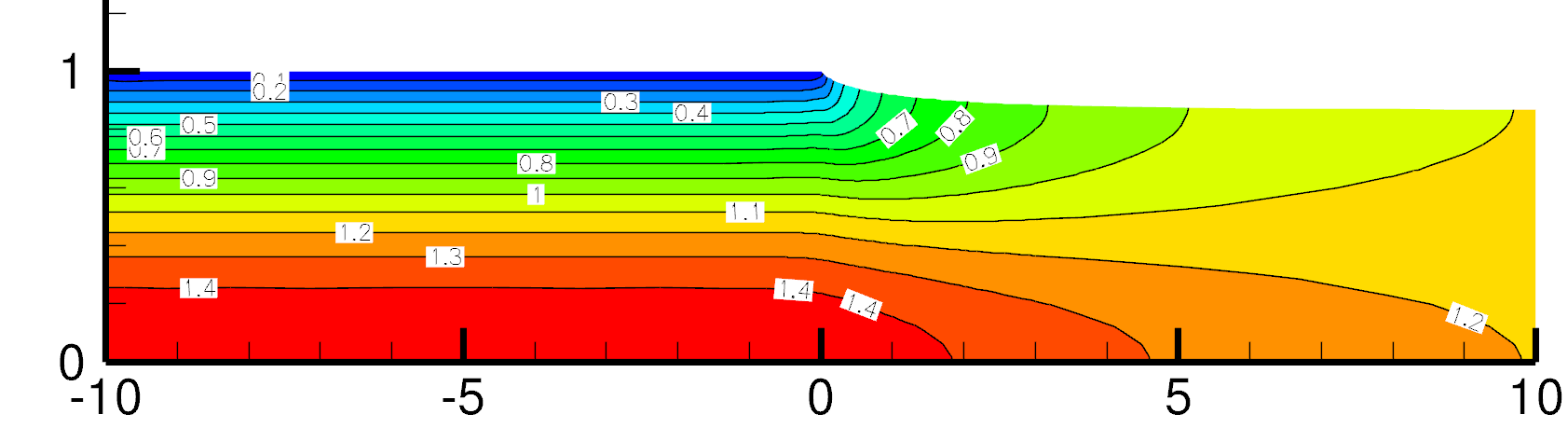}}\\
    & \includegraphics[width=\linewidth]
                    {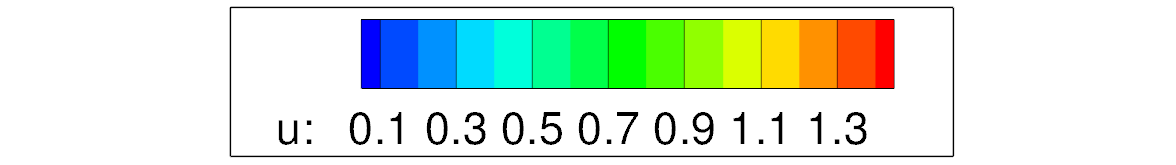}\\
\end{tabular}
\caption{Horizontal velocity component $u$ for $P=10$
for (a) $Re=0$, (b) $Re=3$, (c) $Re=7$, (d) $Re=10$ and (e) $Re=50$. Contours
are indicated at intervals of $0.1$.}
\label{fig:DieSwell contour u}
\end{figure}

\begin{figure}[t]
\centering
\begin{tabular}{ p{0.1\linewidth} p{0.8\linewidth}}
(a) & \raisebox{-0.9\height}{\includegraphics[width=.9\linewidth]
                    {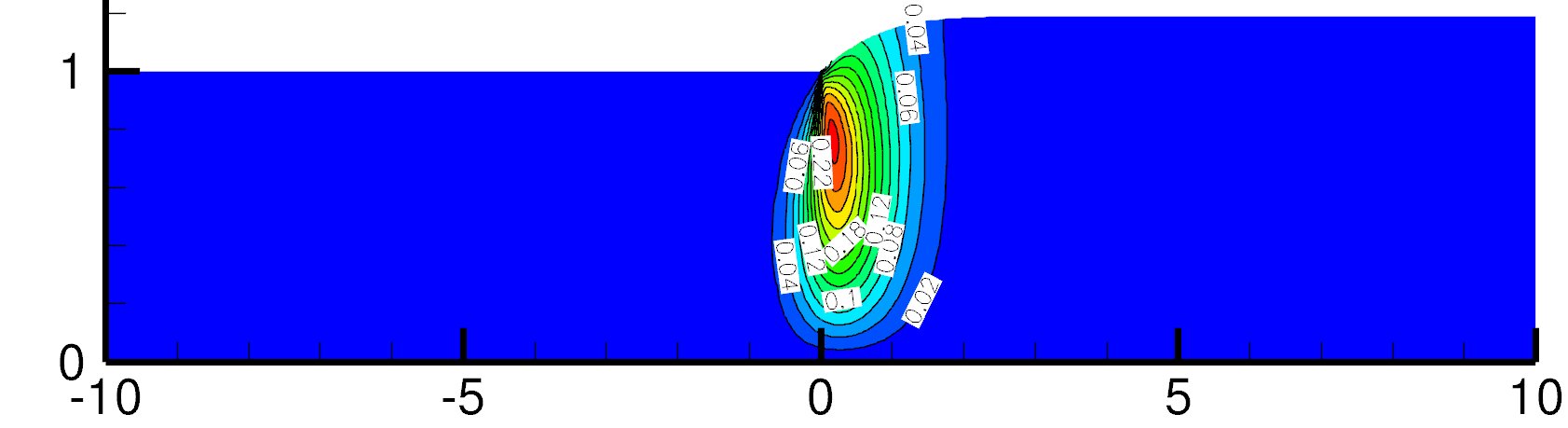}}\\
    & \includegraphics[width=0.5\linewidth]
                    {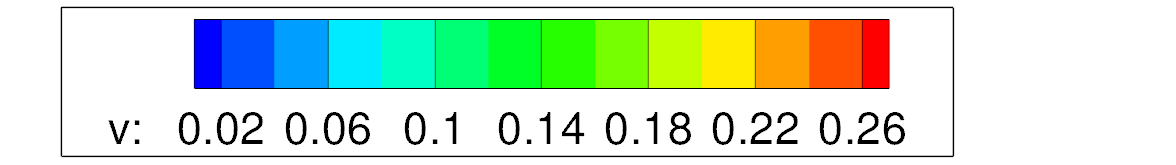}\\
(b) & \raisebox{-0.9\height}{\includegraphics[width=.9\linewidth]
                    {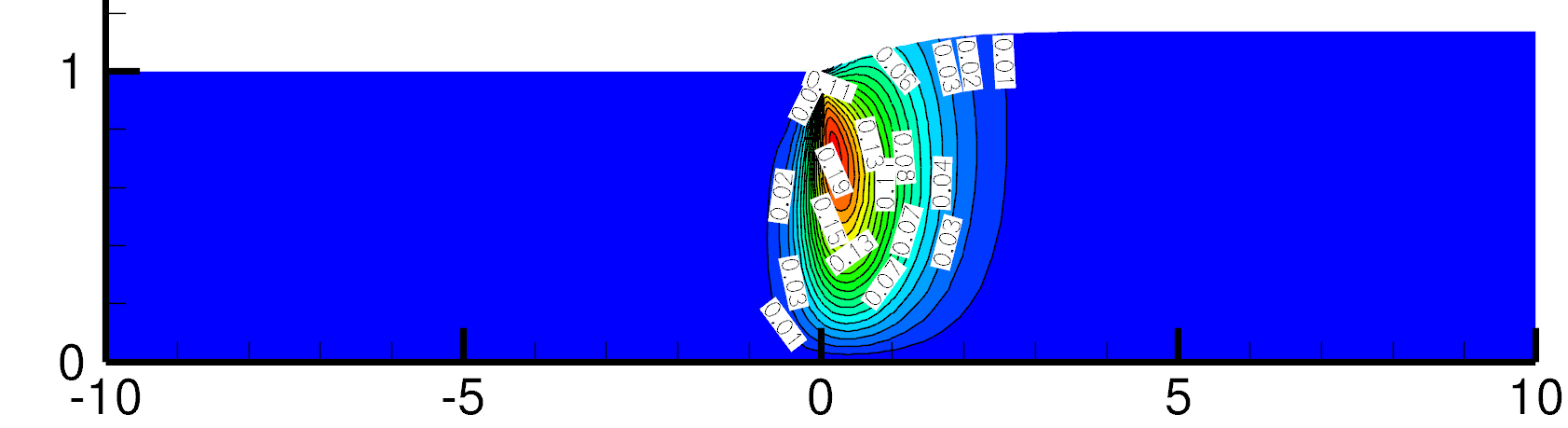}}\\
    & \includegraphics[width=0.5\linewidth]
                    {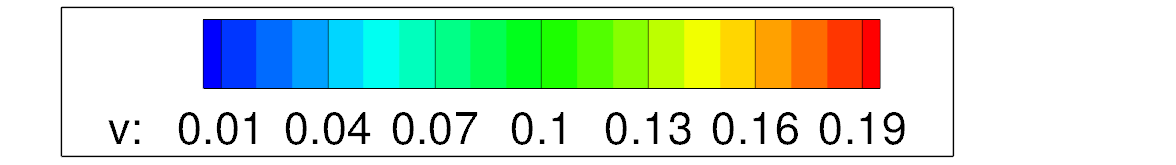}\\
(c) & \raisebox{-0.9\height}{\includegraphics[width=.9\linewidth]
                    {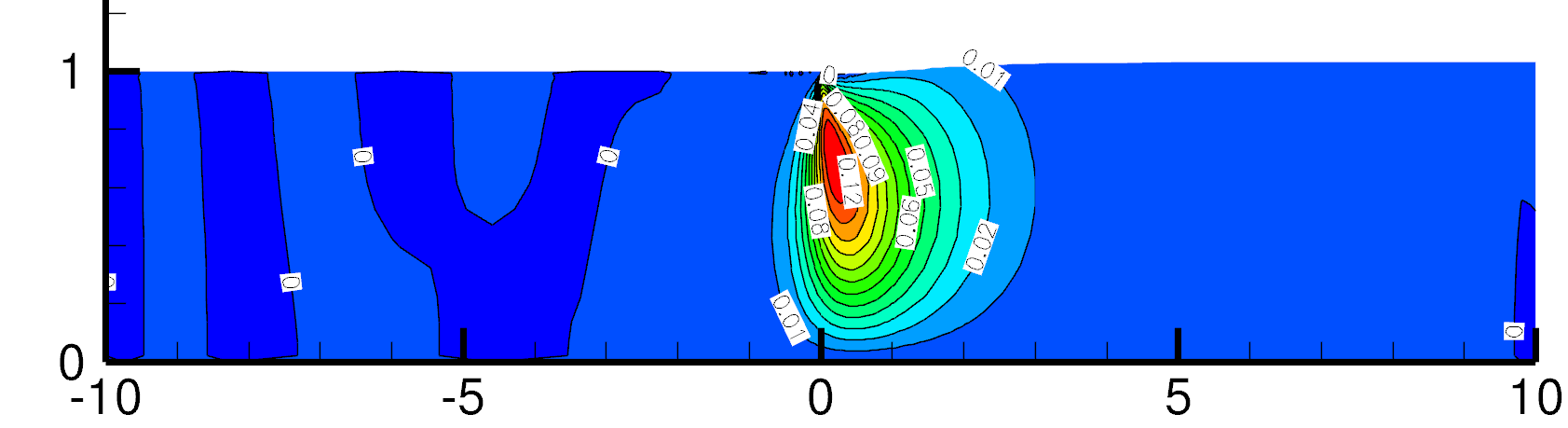}}\\
    & \includegraphics[width=0.5\linewidth]
                    {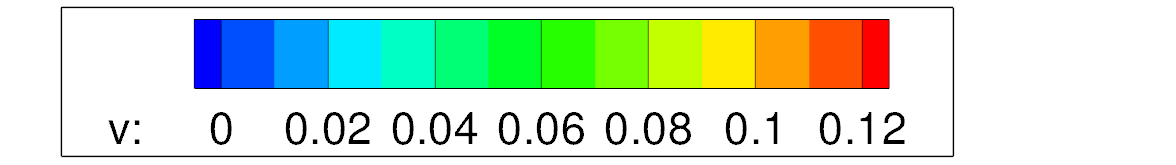}\\
(d) & \raisebox{-0.9\height}{\includegraphics[width=.9\linewidth]
                    {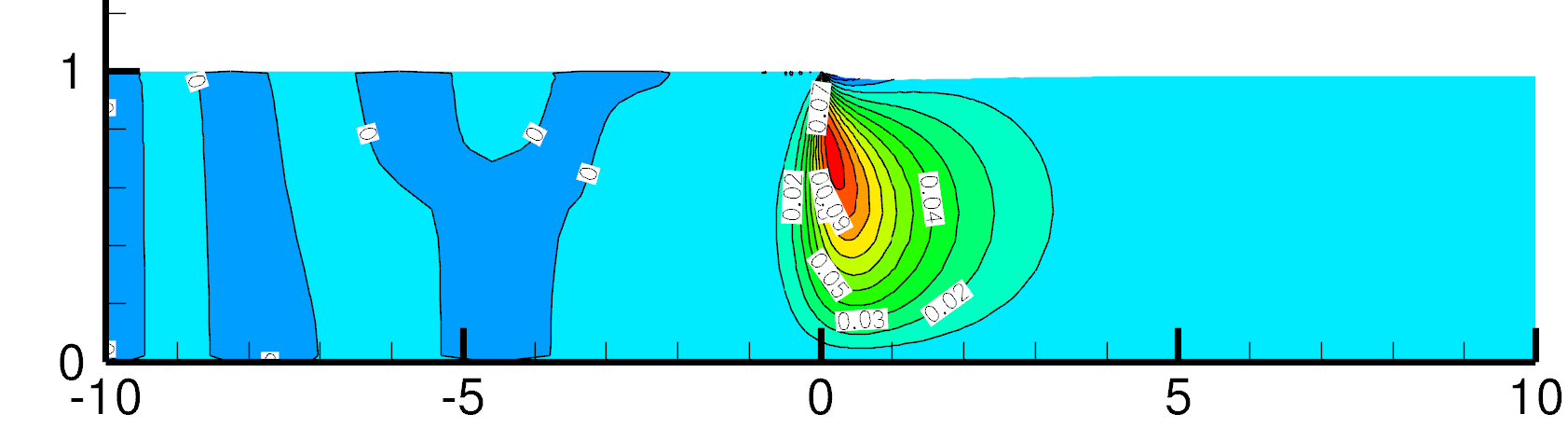}}\\
    & \includegraphics[width=0.5\linewidth]
                    {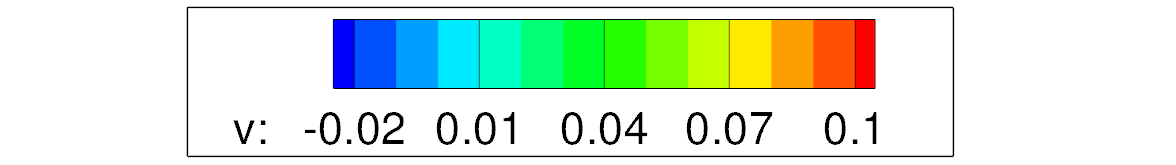}\\
(e) & \raisebox{-0.9\height}{\includegraphics[width=.9\linewidth]
                    {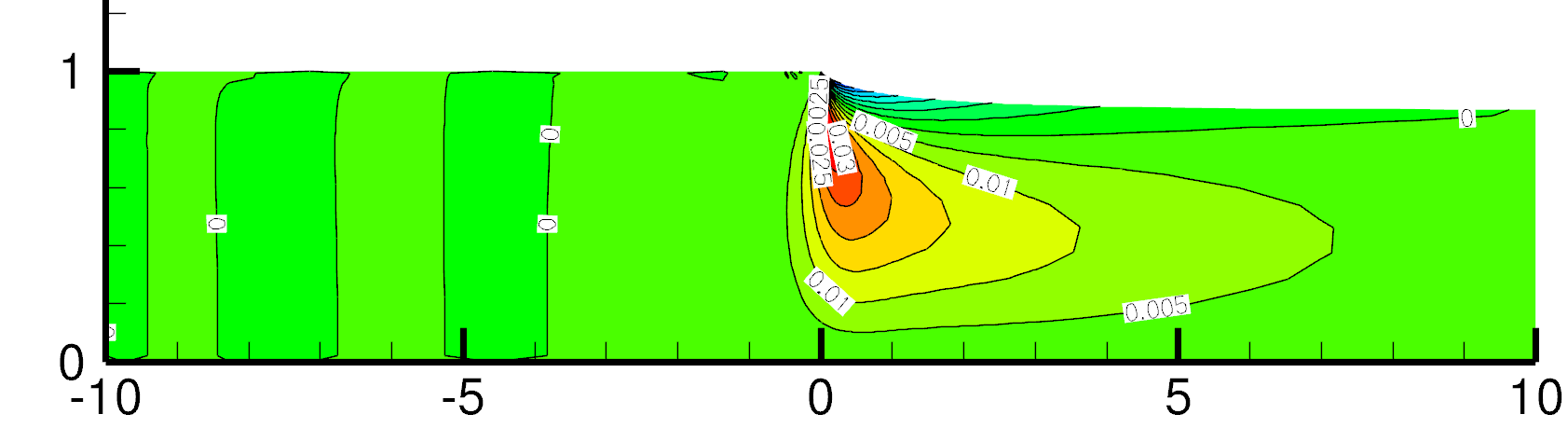}}\\
    & \includegraphics[width=0.5\linewidth]
                    {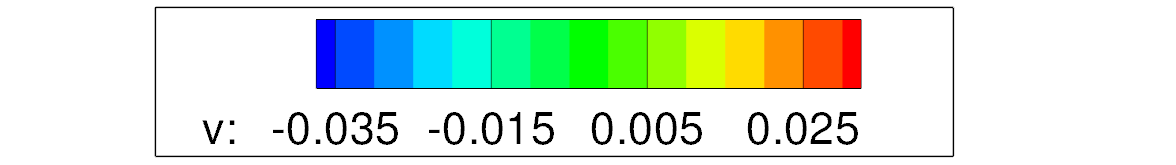}\\
\end{tabular}
\caption{Vertical velocity $v$ for $P=10$
for (a) $Re=0$, (b) $Re=3$, (c) $Re=7$, (d) $Re=10$ and (e) $Re=50$.}
\label{fig:DieSwell v contour}
\end{figure}

We explore the contour plots of the velocity field for a range of Reynolds
numbers in Figures~\ref{fig:DieSwell contour u} (horizontal velocity component
$u$), ~\ref{fig:DieSwell v contour} (vertical velocity component $v$). With
increasing Reynolds number the horizontal velocity increases along the
centreline, the vertical velocity near the singularity induced by the sudden
change in the boundary condition  decreases and the transition zone under the
free surface from Poiseuille flow in the die to plug flow is extended
downstream. This shows that with increasing Reynolds number the particles along
the centreline are accelerated and decelerated near the free surface yielding
the contraction of the free fluid jet. This is indeed the behaviour we would
expect as particles leaving the die will deviate less from their initial path
for increasing inertia.
As pointed out by \cite{Mitsoulis2012a} in order to accommodate the whole
transition zone the domain length of the free fluid jet should be chosen as $L_2
= \Rey$. However, we employ open boundary conditions at outflow which enable us
to compute the extrudate swell accurately in the truncated domain with $L_2=10$.
As demonstrated by \cite{Mitsoulis2011}  the results for extrudate swell with a
domain length $L_2=6$ are virtually identical with those from long domains with
$L_2=\Rey$, for all variables, when using the open boundary condition at
outflow. However, in this case, the swell ratio results are only correct up to
the truncated length as they continuously drop beyond the truncated domain. A
small discrepancy between swell ratios for different domain lengths can
therefore be expected.

\begin{figure}[t]
\centering
\subfloat[Velocity component $u$ along symmetry line
($v=0$).]{\includegraphics[width=.7\linewidth]{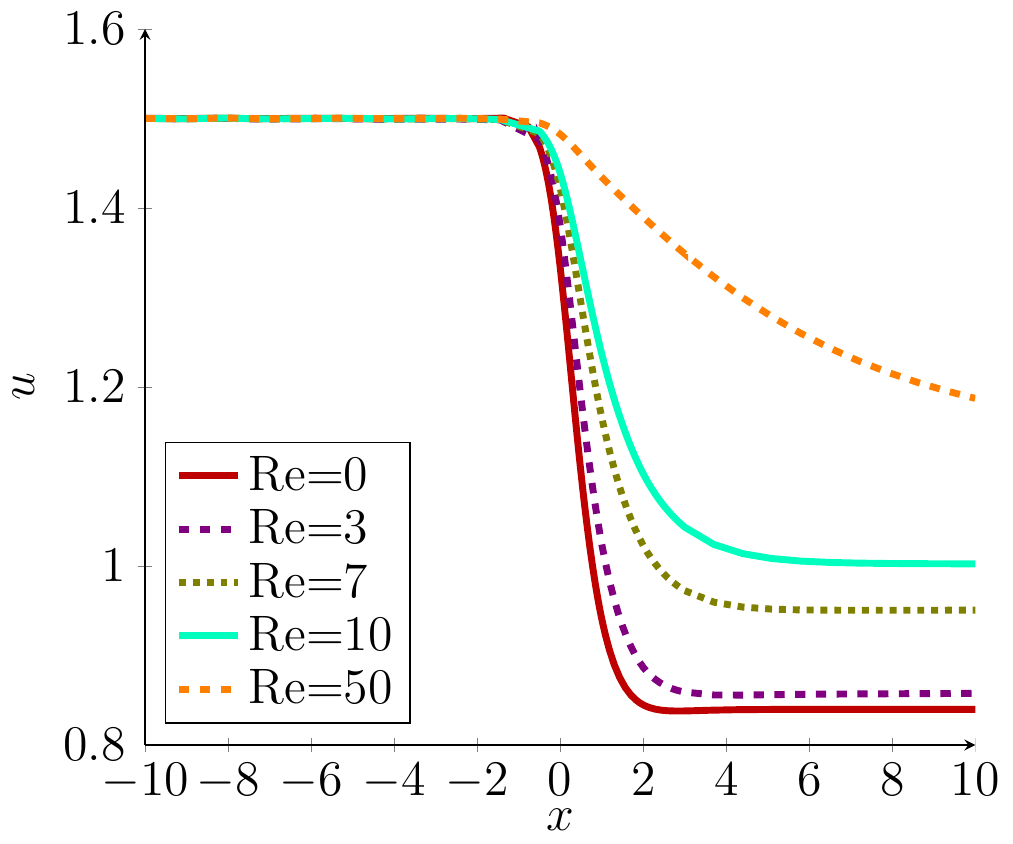}\label{subref:
centreline newt swell u}}\\
\subfloat[Velocity component $u$ along free
surface.]{\includegraphics[width=.5\linewidth]{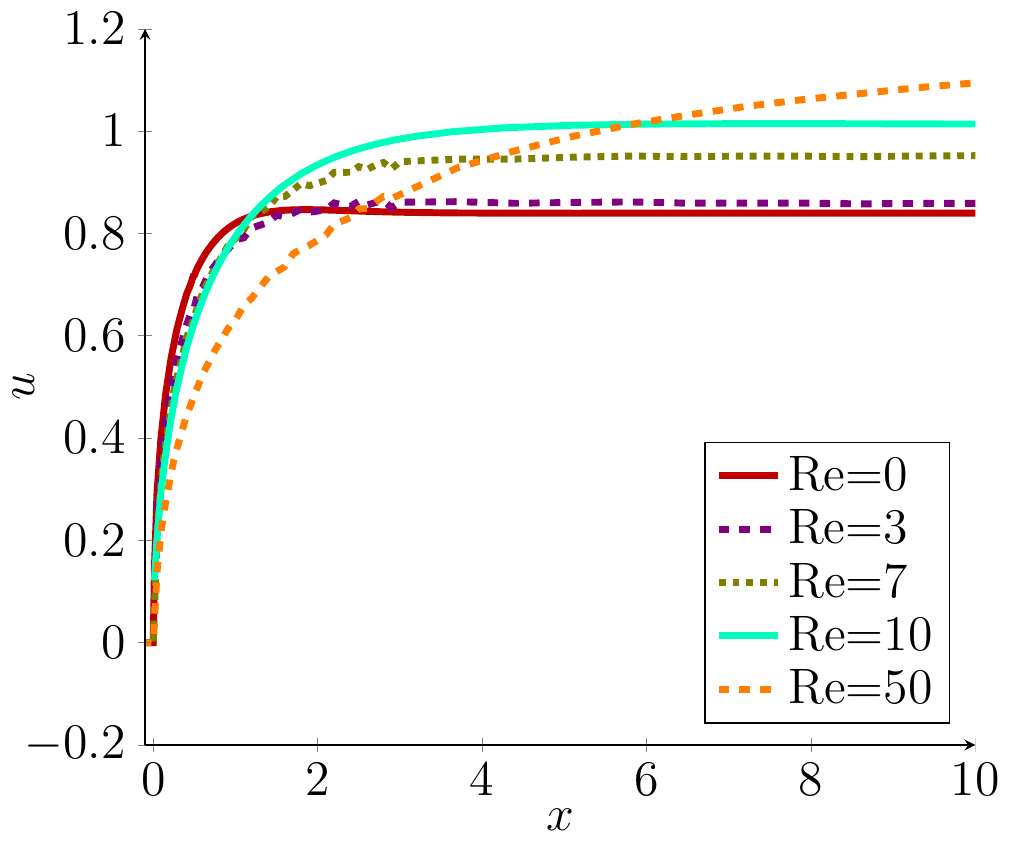}\label{subref:
free surface newt swell u}}
\subfloat[Velocity component $v$ along free
surface.]{\includegraphics[width=.5\linewidth]{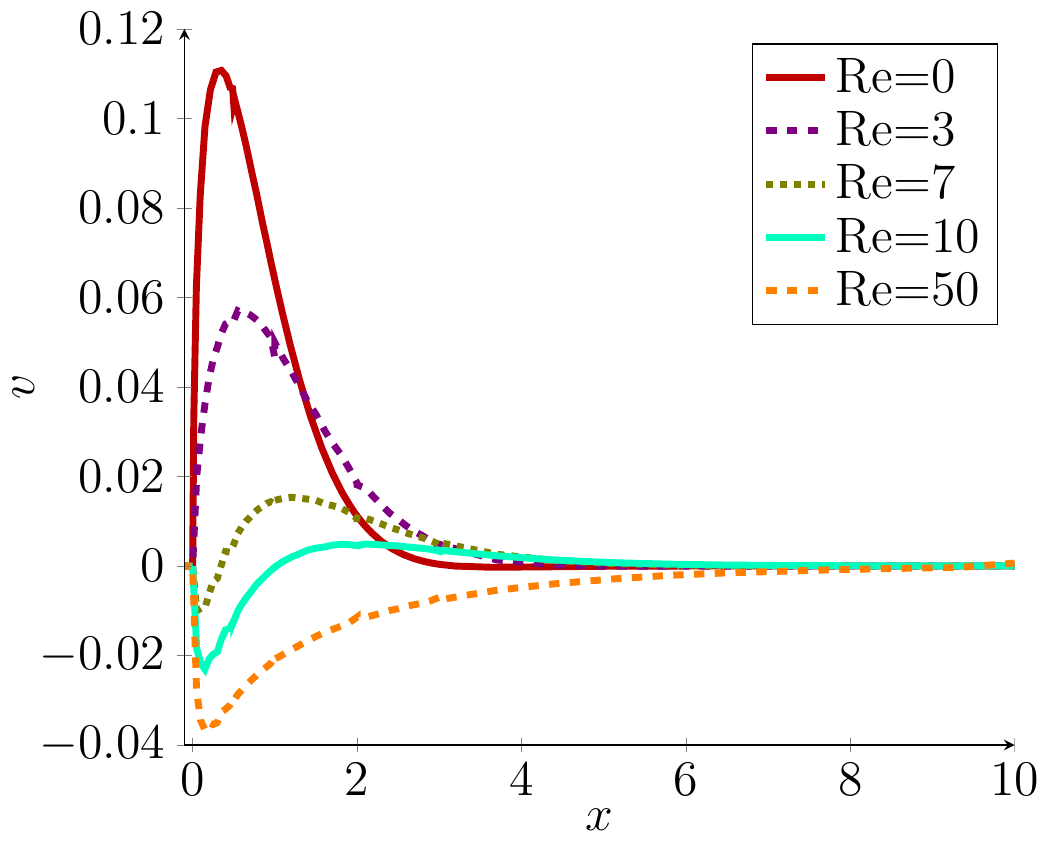}\label{subref:
free surface newt swell v}}
\caption{Dependency of velocity components along \protect\subref{subref: centreline newt swell u} the symmetry line  and \protect\subref{subref: free surface newt swell u}- \protect\subref{subref: free surface newt swell v} along the free surface on the Reynolds number.}
\label{fig:Centreline Newtonian Swell up}
\end{figure}

To investigate the transition from Poiseuille flow to plug flow for increasing
Reynolds number further, we plot the velocity and pressure along different paths
in the domain. Figure~\ref{fig:Centreline Newtonian Swell up} displays the
velocity components along the symmetry line (i.e. $v=0$) and along the free
surface boundary. In Figure~\ref{subref: centreline newt swell u}, we see the
smooth transition of the velocity field from the maximum of the parabolic
profile to the average plug flow velocity given by Equation~\eqref{equ: uplug}, i.e.
$u_{plug}=1/\chi_R$. As the swell decreases with increasing Reynolds number the
plug flow value of the velocity increases with increasing Reynolds number.  With
increasing Reynolds number the change from the maximum parabolic value of the
velocity component $u$ to the plug flow value shifts further downstream. For
$\Rey=0$, the velocity reaches the plug flow value at around $x\approx 3$, for
$\Rey=10$ at $x \approx 6$ and for $\Rey=50$ the plug flow value is not reached
within our computational domain. However, as pointed out above, due to the use
of open boundary conditions at outflow, the velocity and pressure profiles stay
accurate even if they are truncated at outflow. 

Along the free surface boundary
(Figure~\ref{subref: free surface newt swell u}), the velocity component $u$ increases sharply near the
die exit until it reaches the plug flow value while the velocity component $v$ goes
through a maximum near the die exit for $\Rey=0$ and $\Rey=3$ and through a
minimum for $\Rey>7$, when particles are no longer constrained by the no-slip
boundary condition (Figure~\ref{subref: free
surface newt swell v}). This causes the swell (for $v>0$) or the contraction (for
$v<0$) of the free surface near the die exit until the surface is sufficiently
curved to obtain a zero total shear stress (i.e. $\mathbf{t} \cdot \cauchystress
\cdot \mathbf{n}=0$). Further downstream when the free surface boundary has
reached its maximum swelling value, the vertical velocity component reaches zero
in accordance with the condition of no particle penetration along the surface
(horizontal free surface boundary has outward normal $\mathbf{n}=(0,1)$ and therefore
$\mathbf{u} \cdot \mathbf{n} = v =0$). The maximum value of $v$ along the free
surface decreases with increasing Reynolds number ($0\leq \Rey \leq 5$). For the
range of Reynolds number that causes a delayed die swell the velocity component
$v$ first undergoes a sharp minimum and then goes through a maximum ($6\leq \Rey
\leq 10$). For the range of Reynolds numbers that cause a contraction of the
free Newtonian jet, the velocity component $v$ goes through a minimum and then
slowly approaches zero ($\Rey > 10$). 

\begin{figure}[t]
\centering
\subfloat[Velocity component $u$ in cross stream
direction.]{\includegraphics[width=.7\linewidth]{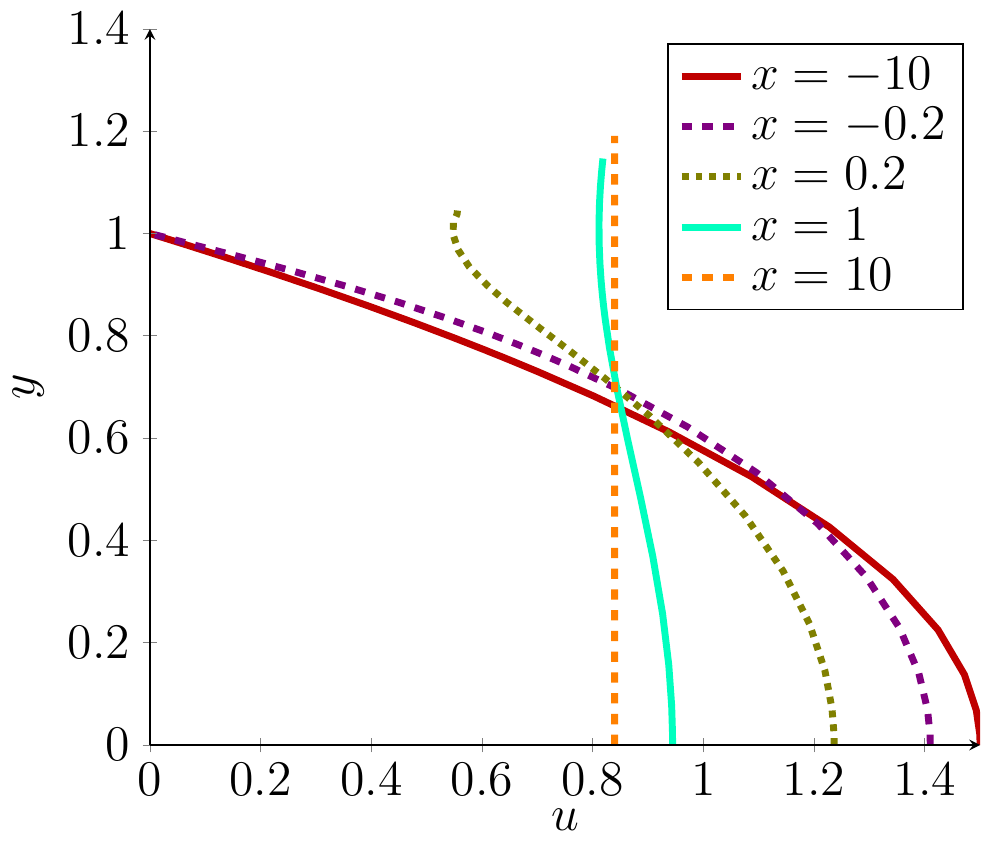}}\\
\subfloat[Velocity component $v$ in cross stream
direction.]{\includegraphics[width=.7\linewidth]{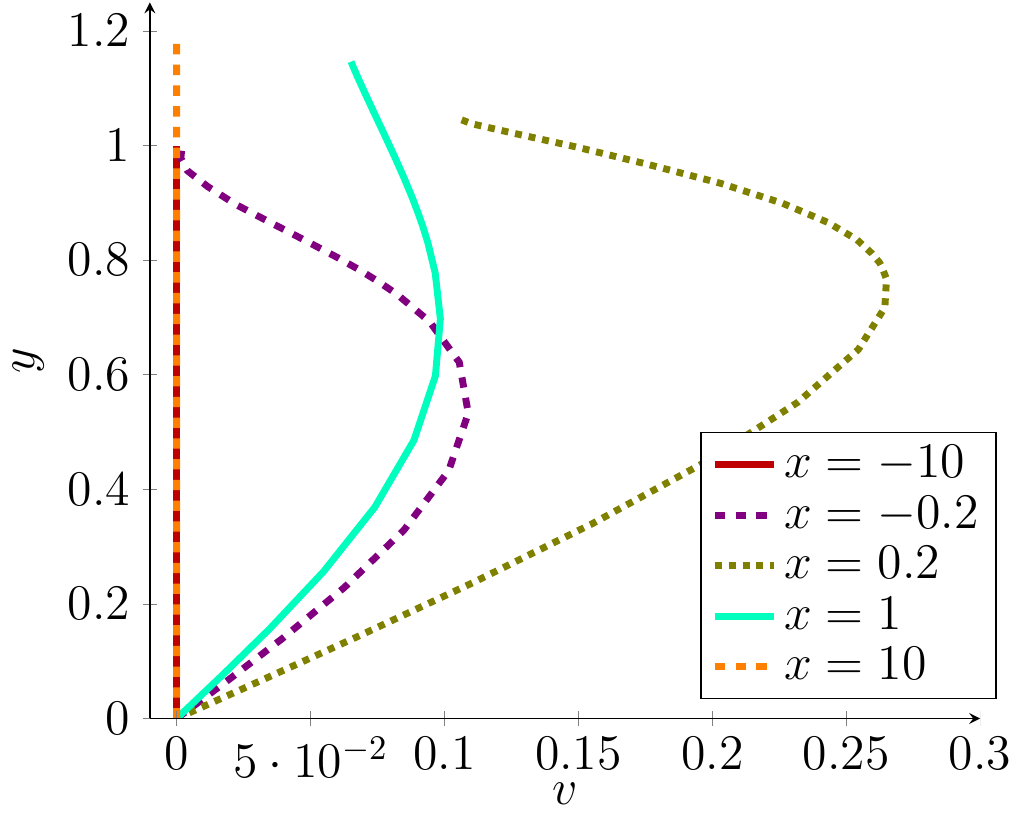}}\\
\caption{Velocity components in cross stream direction at inflow ($x=-10$), near the die exit ($x=-0.2$, $x=-0.2$), further downstream in the free jet region $x=1$ and at outflow $x=10$.}
\label{fig: swell Newton cross streamwise}
\end{figure}

Figure~\ref{fig: swell Newton cross streamwise} shows the velocity components in
the cross stream wise direction at inflow ($x=-10$), near the die exit
($x=-0.2$, $x=0.2$), further downstream in the free jet region $x=1$ and at
outflow $x=10$. The velocity component $u$, is parabolic at inflow, shortly
before the die exit ($x=-0.2$) the parabolic profile flattens inside the die,
after the die exit the parabolic profile flattens further and builds a boundary
layer in which it goes through a minimum $x=0.2$, then flattens increasingly
until the plug flow value is reached. The vertical velocity component, which is
zero at inflow, forms a parabolic-like profile with a small boundary layer near
the die exit inside the die, which first sharpens shortly after exiting the die
and then relaxes back to the zero value.

\begin{figure}[t]
\centering
\begin{tabular}{ p{0.1\linewidth} p{0.8\linewidth}}
(a) & \raisebox{-0.9\height}{\includegraphics[width=\linewidth]
                    {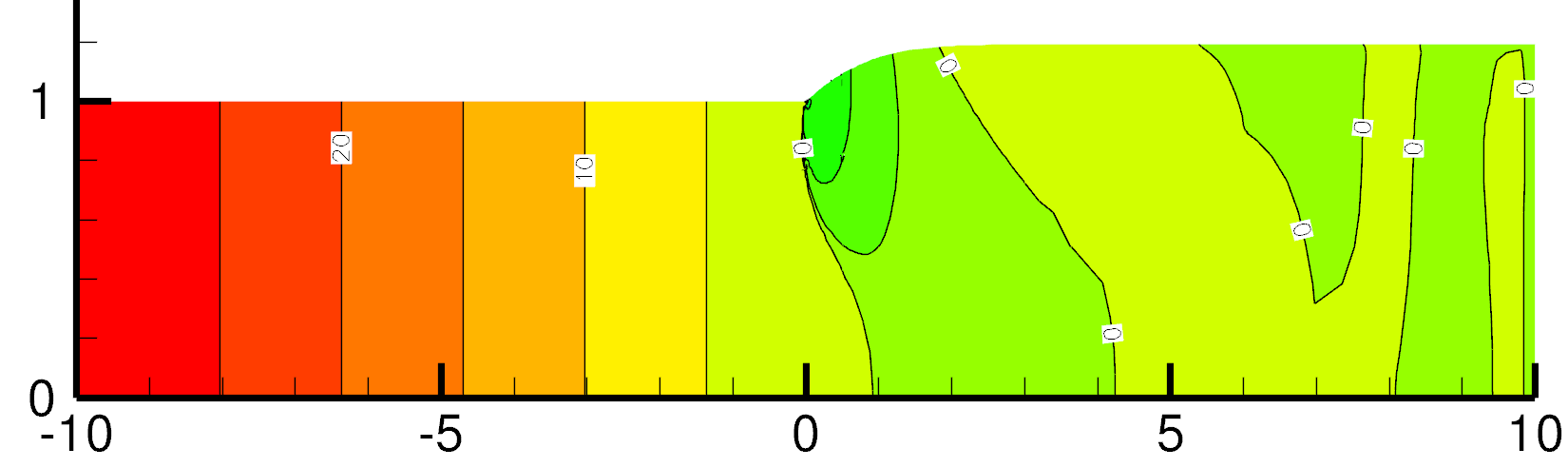}}\\
(b) & \raisebox{-0.9\height}{\includegraphics[width=\linewidth]
                    {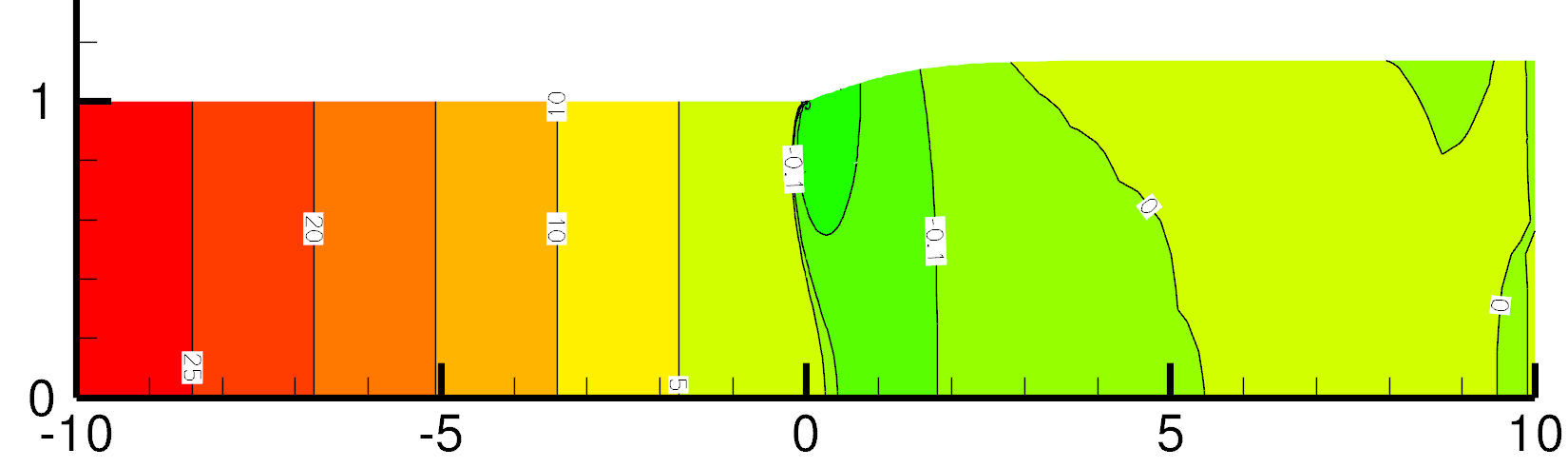}}\\
(c) & \raisebox{-0.9\height}{\includegraphics[width=\linewidth]
                    {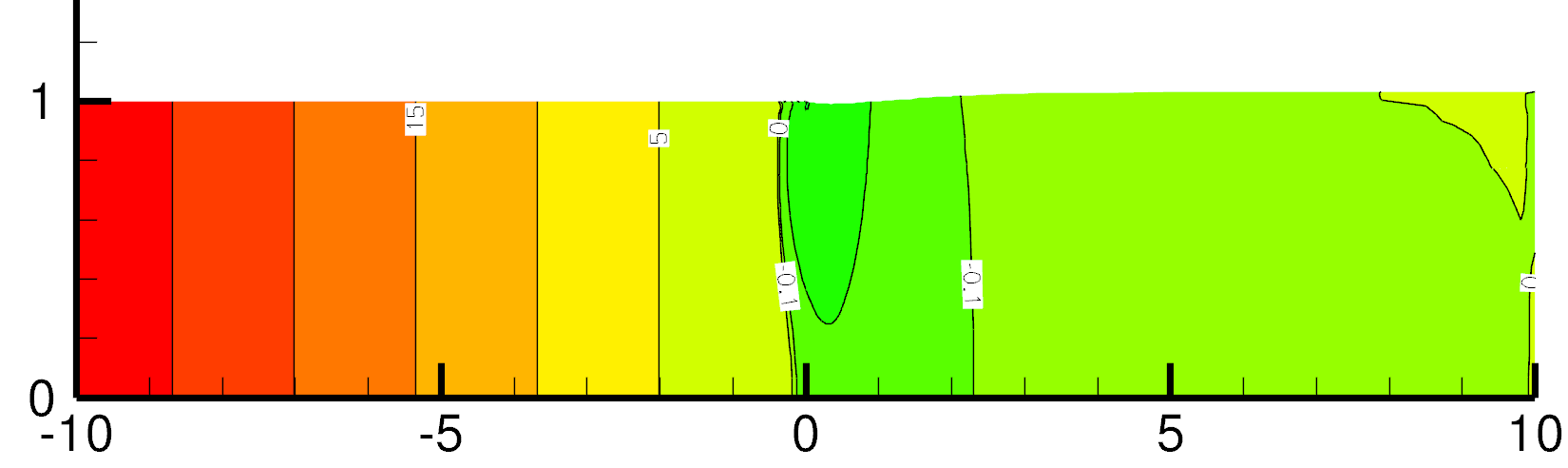}}\\
(d) & \raisebox{-0.9\height}{\includegraphics[width=\linewidth]
                    {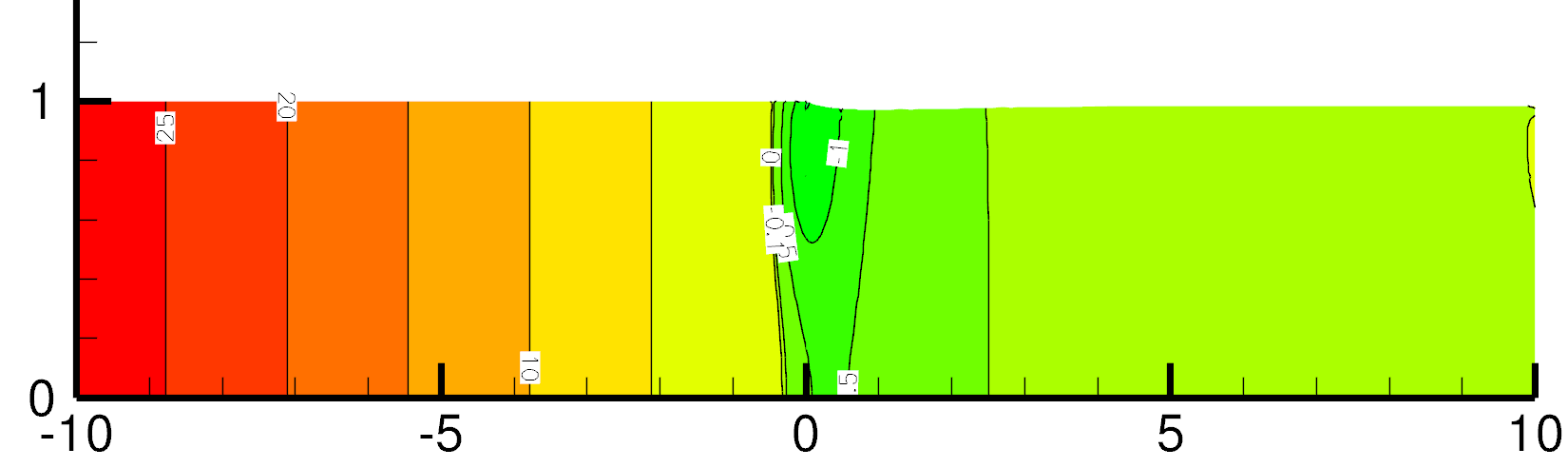}}\\
(e) & \raisebox{-0.9\height}{\includegraphics[width=\linewidth]
                    {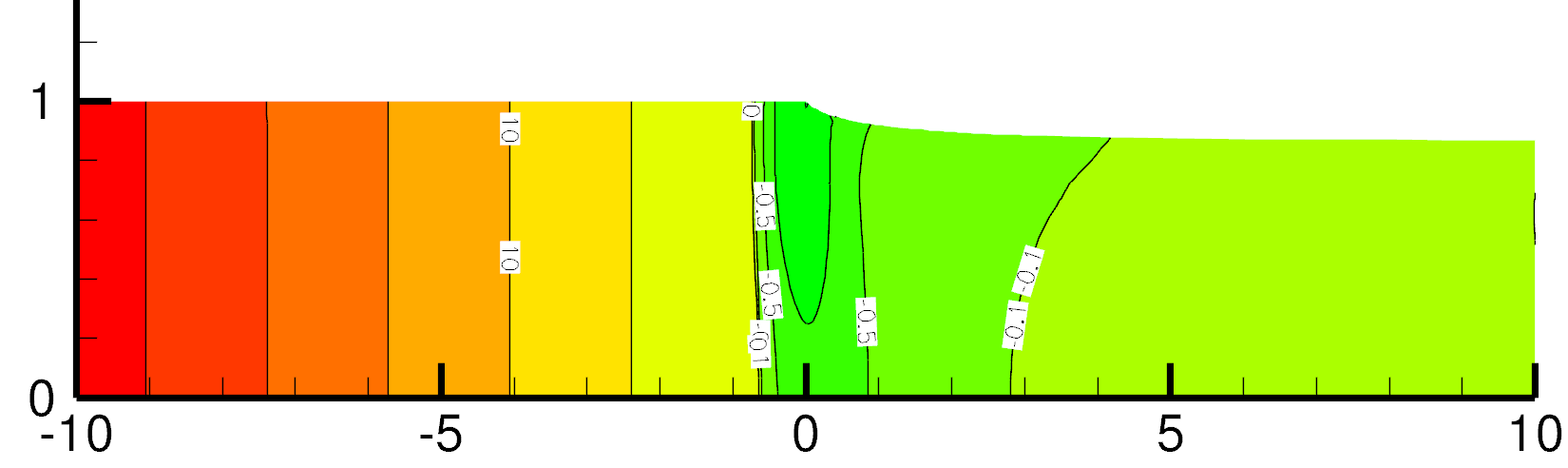}}\\
    & \includegraphics[width=0.5\linewidth]
                    {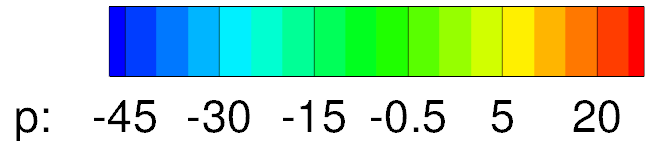}\\
\end{tabular}
\caption{Pressure $p$ for $P=10$
for (a) $Re=0$, (b) $Re=3$, (c) $Re=7$, (d) $Re=10$ and (e) $Re=50$. Contours
are indicated at intervals of $0.1$.}
\label{fig:DieSwell p contour}
\end{figure}

\begin{figure}[t]
\centering
\subfloat[Pressure $p$ along
centreline.]{\includegraphics[width=.7\linewidth]{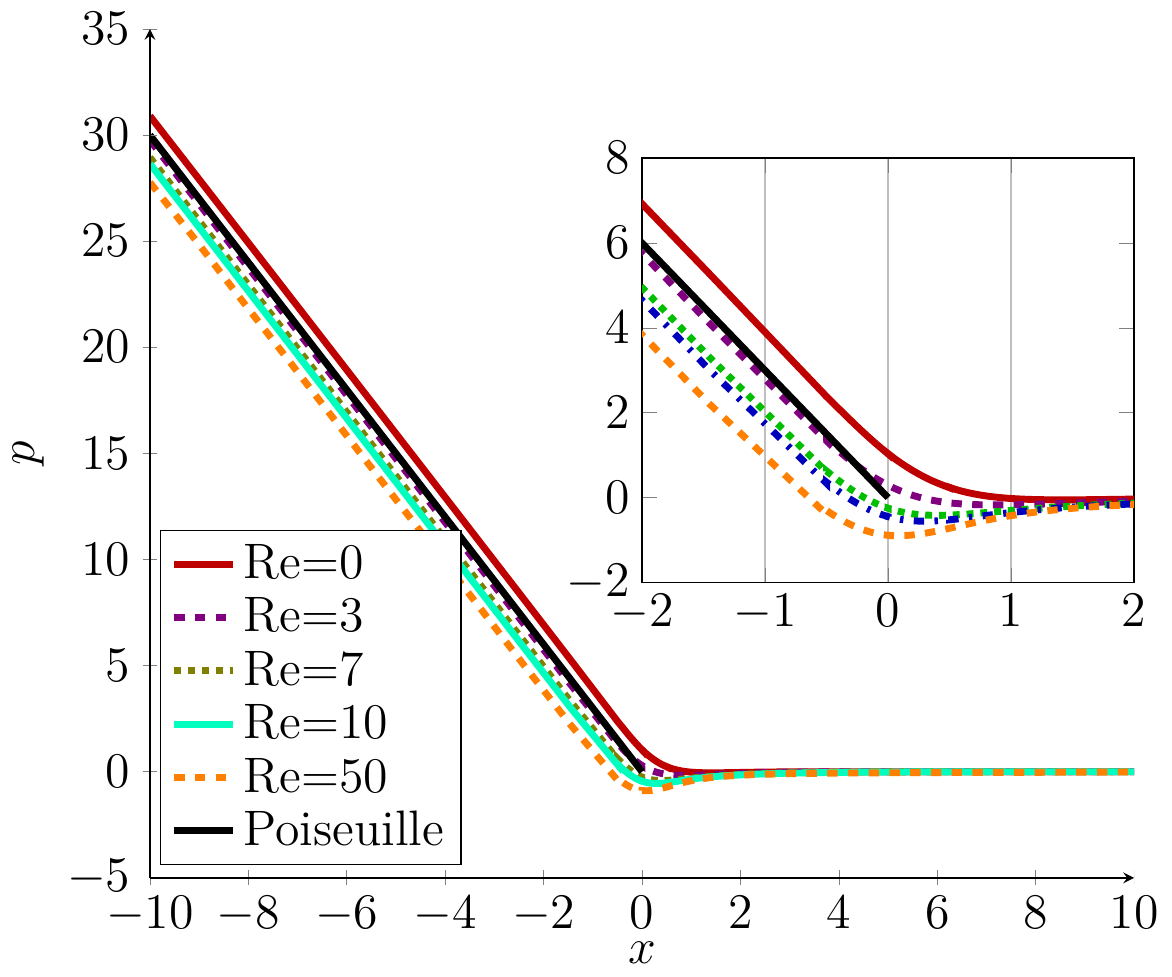}\label{subref:
centreline newt swell p}}\\
\subfloat[Pressure $p$ around the singularity at the die
exit.]{\includegraphics[width=.7\linewidth]{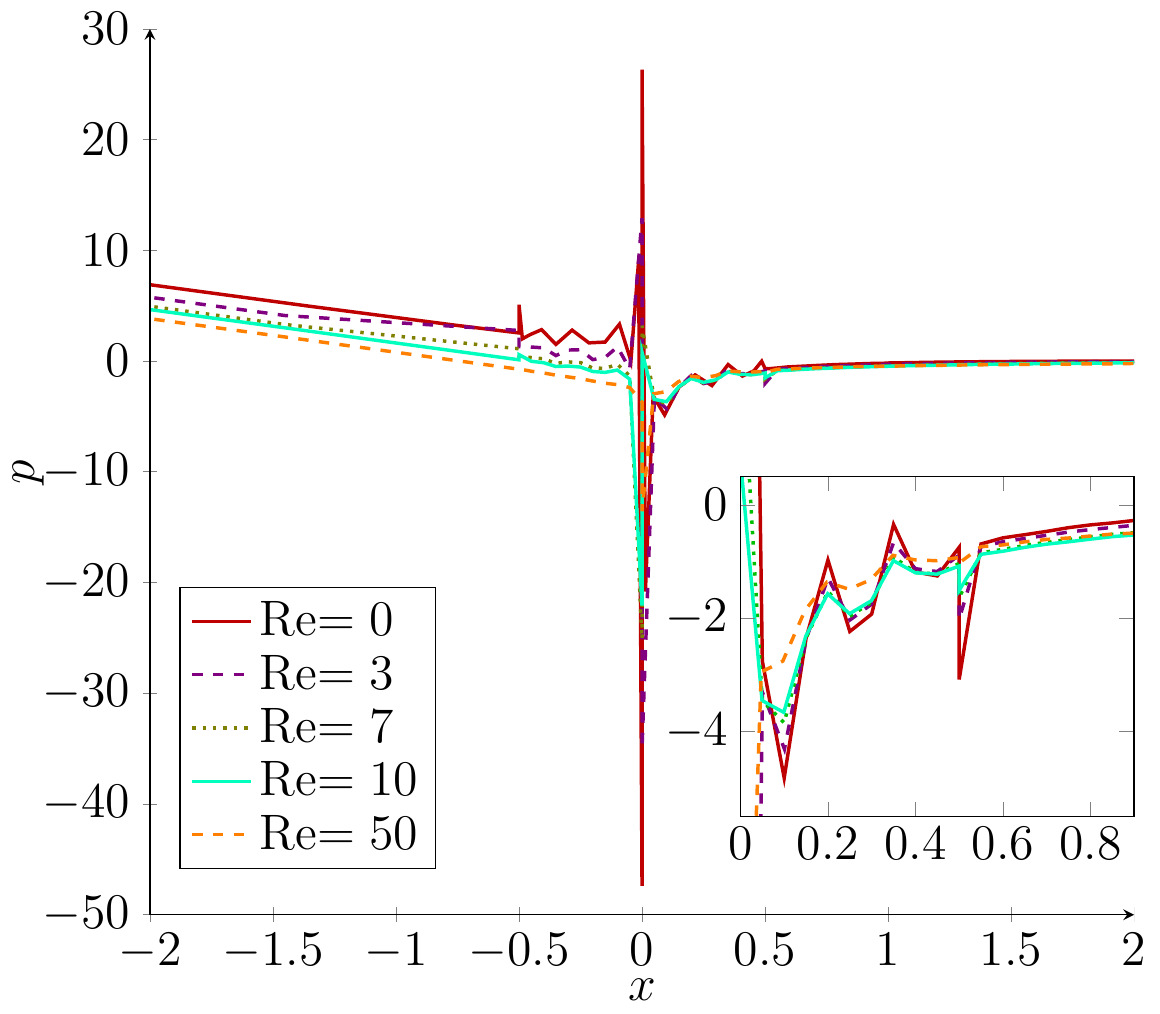}\label{subref:
fs newt swell p Re}}\\
\caption{Plots of pressure $p$ along \protect\subref{subref: centreline newt swell p} the centreline and \protect\subref{subref: fs newt swell p Re} the wall and the free surface.}
\label{fig:pathplots newt swell p Re}
\end{figure}

In the contour plots for the pressure $p$ displayed in 
Figure~\ref{fig:DieSwell p contour},  we observe that the pressure isobars are
curved near the die exit and in the free jet region into the downstream direction for low Reynolds number
($\Rey=0,3,7$) and into the upstream direction for higher Reynolds numbers
($\Rey>10$). The change in the pressure becomes more apparent when we explore
the pressure values along the symmetry line 
(Figure~\ref{subref: centreline newt swell p}). Inside the die, the pressure
gradient is constant as expected for Poiseuille flow. However, near the die exit ($x=0$) the pressure smoothly approaches zero for the plug flow. For higher Reynolds numbers the pressure
on the centreline goes through a minimum.
This behaviour of the pressure yields a shift in the pressure values at inflow,
which is expressed by the pressure exit correction as defined in
Equation~\eqref{equ: pressure exit correction def}.

\subsection{Impact of slip}

\begin{figure}[t]
\centering
\subfloat[]{\includegraphics[width=.45\linewidth]{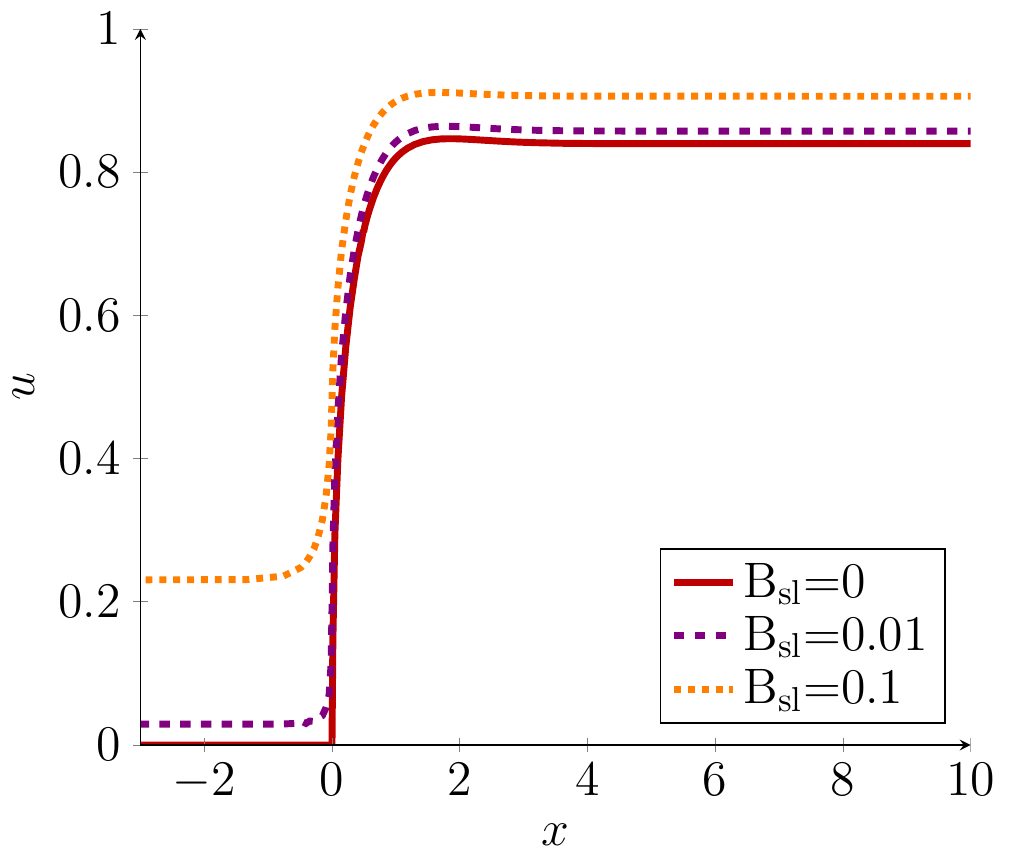}\label{subfig:
newton swell slip fs u}}
\subfloat[]{\includegraphics[width=.45\linewidth]{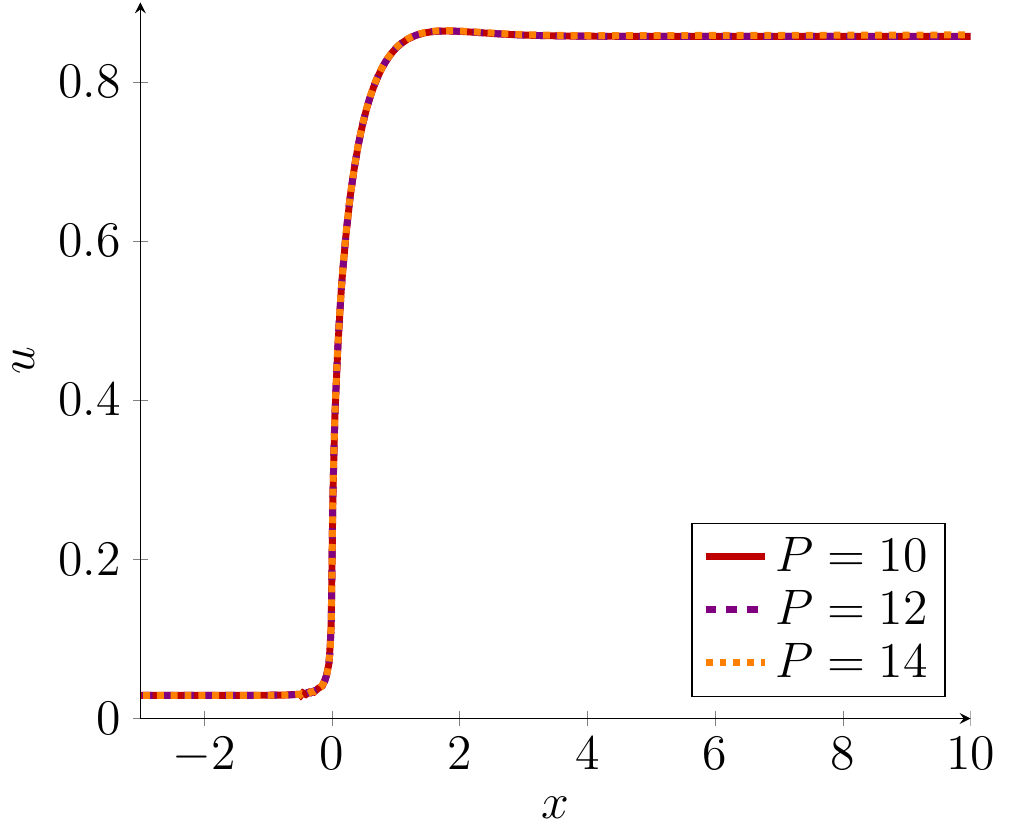}\label{subfig:
newton swell slip fs u P}} \\
\subfloat[]{\includegraphics[width=.45\linewidth]{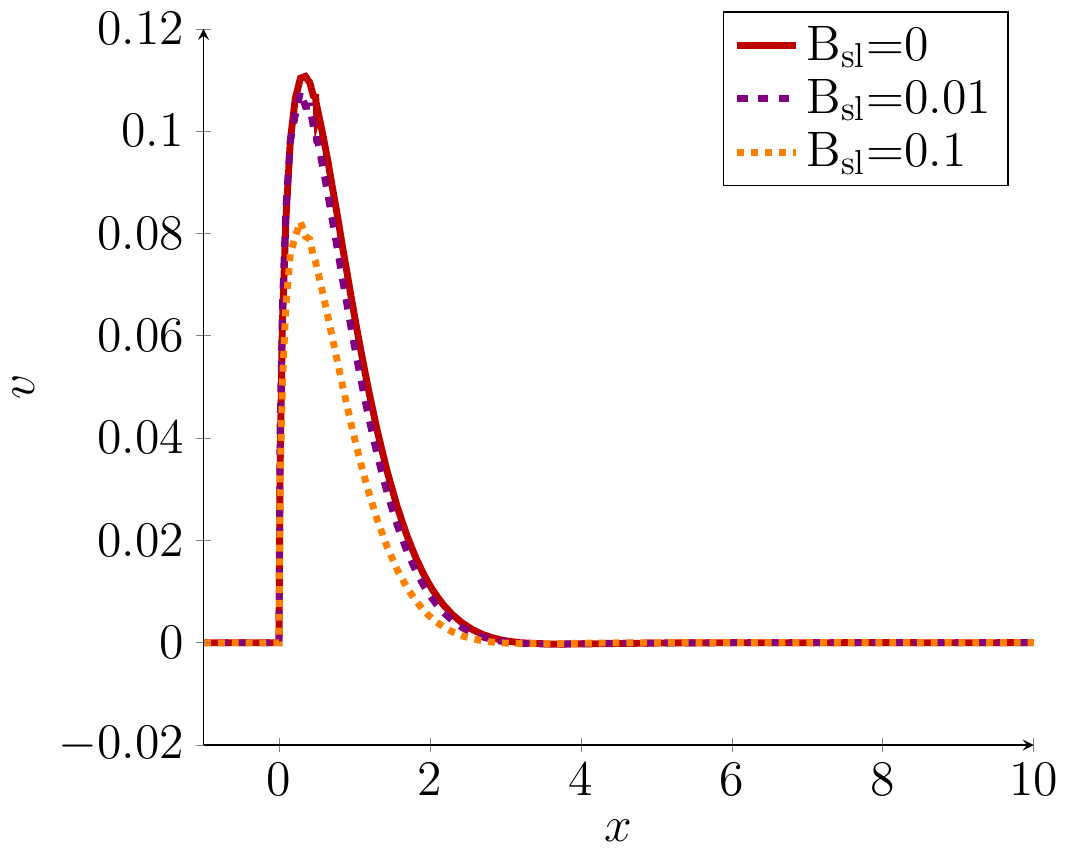}\label{subfig:
newton swell slip fs v}}
\subfloat[]{\includegraphics[width=.45\linewidth]{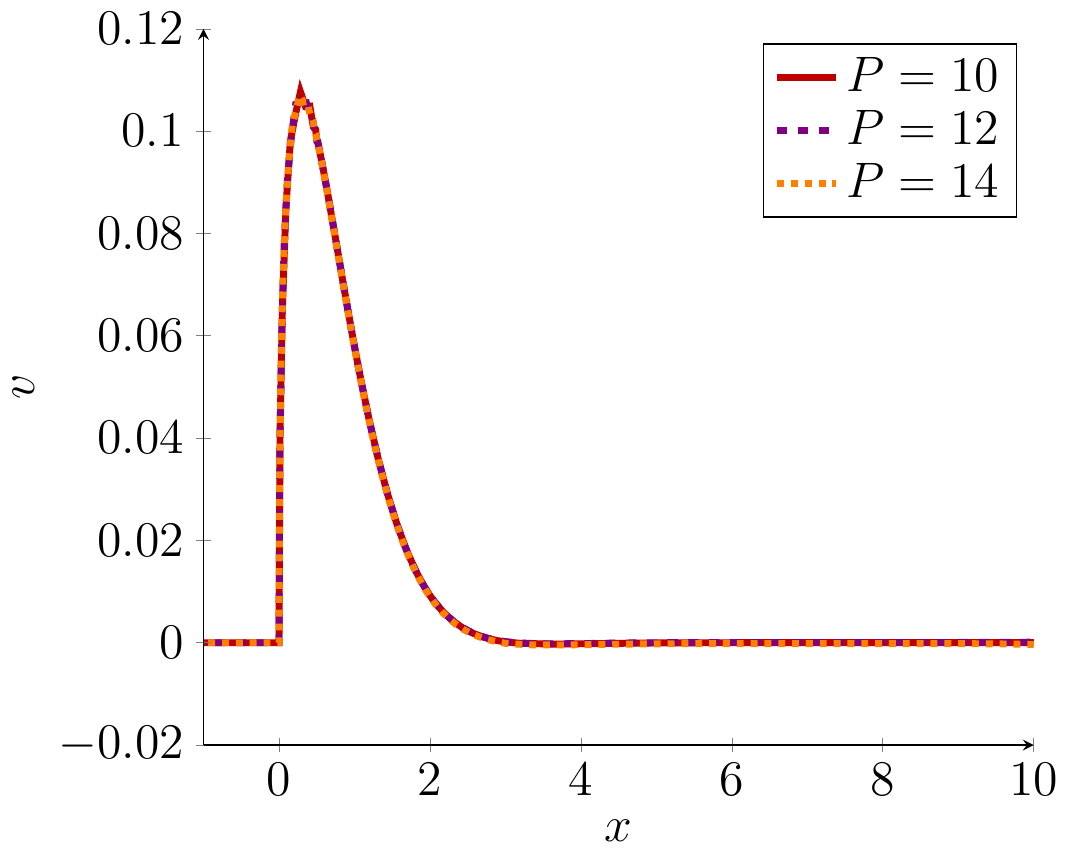}\label{subfig:
newton swell slip fs v P}}\\
\subfloat[]{\includegraphics[width=.45\linewidth]{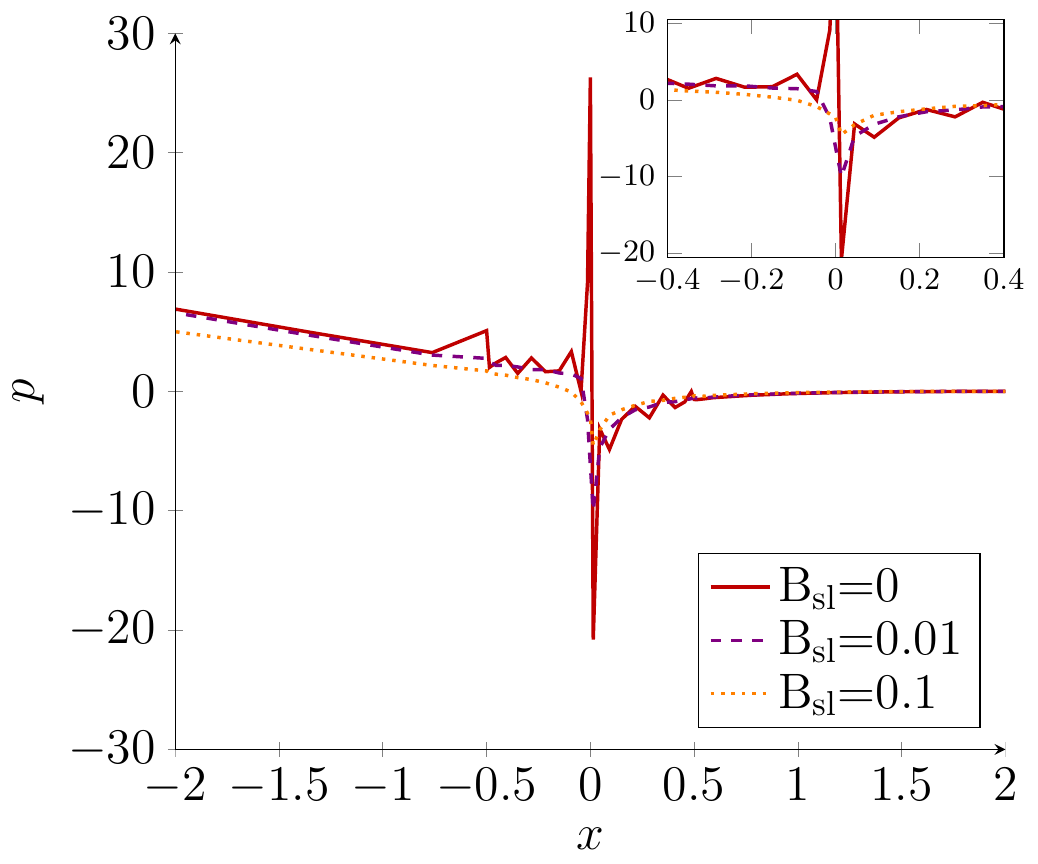}\label{subfig:
newton swell slip fs p}}
\subfloat[]{\includegraphics[width=.45\linewidth]{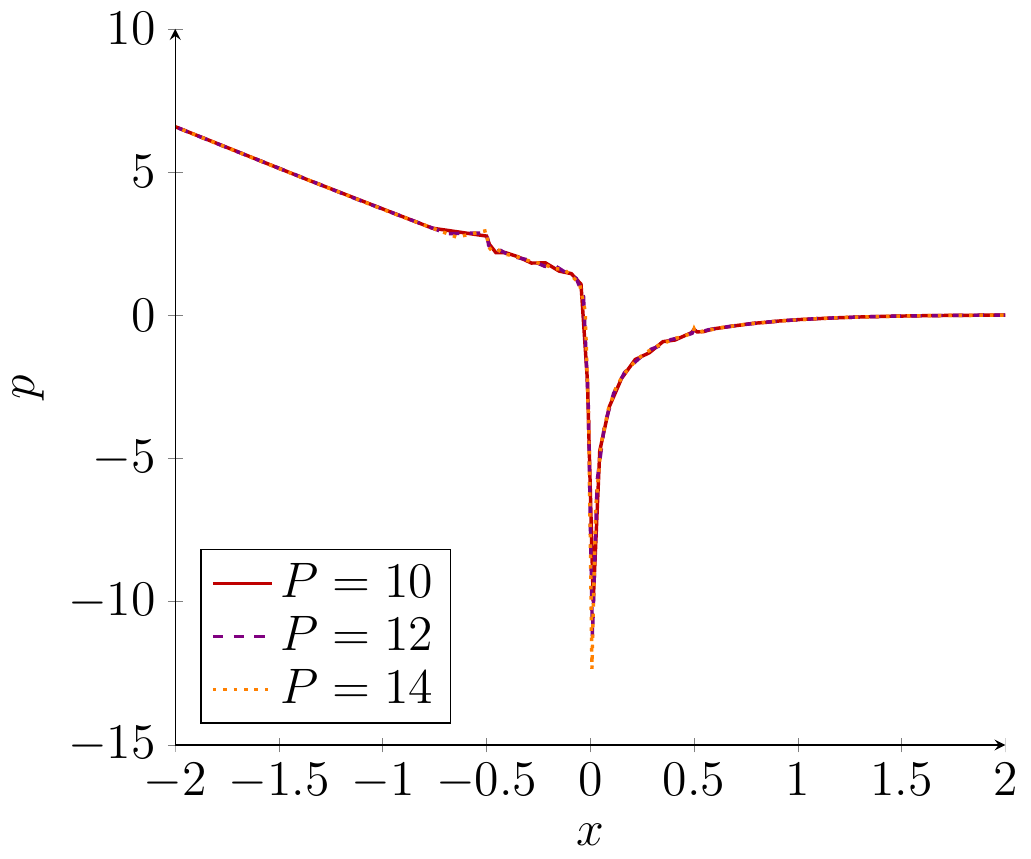}\label{subfig:
newton swell slip fs p P}}\\
\caption{Dependence of \protect\subref{subfig: newton swell slip fs u} velocity components $u$, \protect\subref{subfig: newton swell slip fs v} $v$ and \protect\subref{subfig: newton swell slip fs p} pressure on the slip parameter for $P=10$ and on mesh refinement for $B_{sl}=0.01$ (\protect\subref{subfig: newton swell slip fs u P},\protect\subref{subfig: newton swell slip fs v P},\protect\subref{subfig: newton swell slip fs p P}). }
\label{fig: newton swell slip fs}
\end{figure}

To alleviate the pressure singularity at the die exit, we investigate the effect
of slip along the die wall on the dependent variables for $\Rey=0$. We therefore change the
inflow profile according to Equation~\eqref{equ: slip inflow} and employ the
slip condition \eqref{equ: slip boundary} along the die wall. We explore the
velocity field and the pressure along the free surface for slip parameter values  of
$B_{sl}=0.01$, $B_{sl}=0.1$ and $B_{sl}=0$ (no-slip) in Figure~\ref{fig: newton
swell slip fs}. With the introduction of slip along the wall, the horizontal
velocity component experiences a smooth transition at the die exit in vast
contrast to the kink at the singularity that is observed for the no-slip
condition ($B_{sl}=0$) along the wall (Figure~\ref{subfig: newton swell slip fs
u}). The change for the vertical velocity remains sudden and features a kink at
the singularity. However, the maximum value of the vertical velocity component
decreases with increasing slip (Figure~\ref{subfig: newton swell slip fs v}).

   \begin{table}[t]
      \centering
      \ra{1.2}
      \footnotesize
   \caption{Dependence of the swelling ratio on $P$ for $B_{sl}=0.1$ and $B_{sl}=0.01$.}
   \rowcolors{2}{gray!25}{white}
   \begin{tabular}{rrr}
   \mytoprule
   $P$ & $B_{sl}=0.1$ & $B_{sl}=0.01$ \\ 
   \mymidrule
   10 & 1.1041 & 1.1671 \\ 
   12 & 1.1041 & 1.1673 \\ 
   14 & 1.1040 & 1.1670 \\ 
   \cite{Mitsoulis2012a} & 1.1041 & 1.1708 \\
   \mybottomrule
   \end{tabular}
   \label{tab: swelling ratio slip P}
   \end{table}

\begin{figure}[t]
\centering
\includegraphics[width=.8\linewidth]{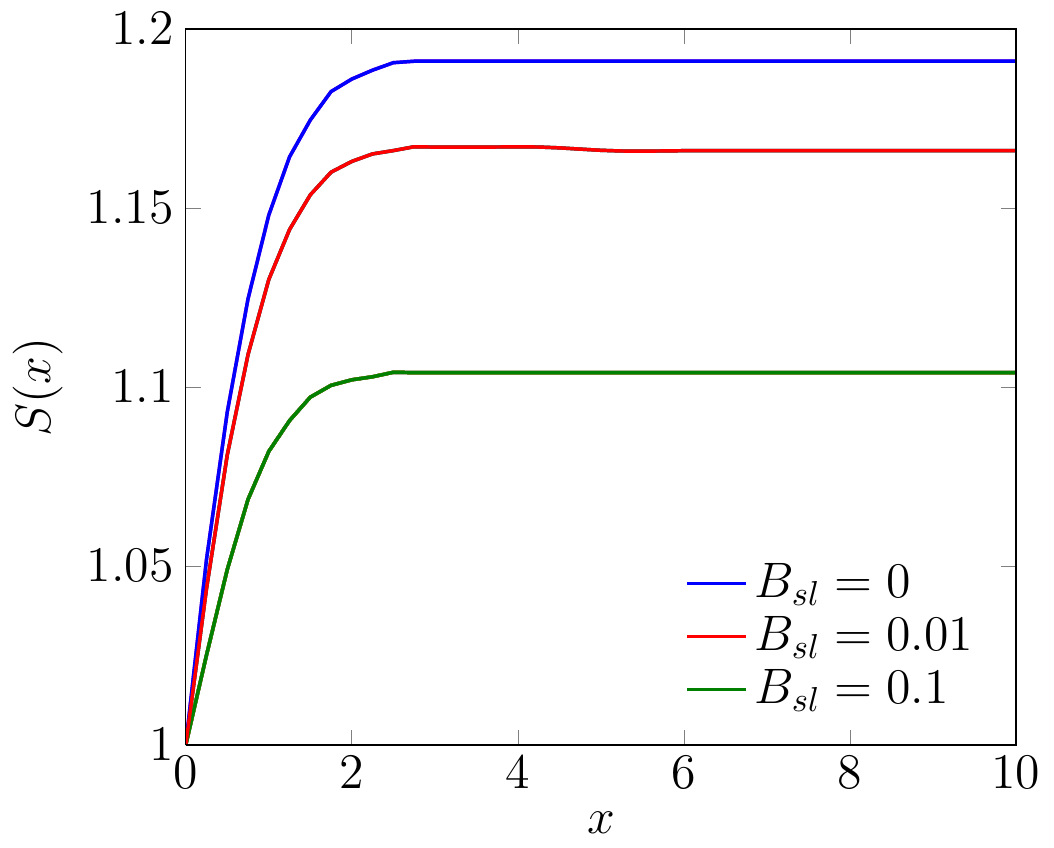}
\caption{Free surface spline profile for increasing slip parameter.}
\label{fig:newt swell spline slip}
\end{figure}

The pressure profile at the singularity is changed drastically with slip along
the wall and the Gibbs oscillations disappear (Figure~\ref{subfig: newton swell
slip fs p}, \protect\subref{subfig: newton swell slip fs p P}). Even though the
minimum of the pressure does not show a converging trend in the range of the
employed polynomial orders, its value only increases slightly with increasing
$P$ (Figure~\ref{subfig: newton swell slip fs p P}). Table~\ref{tab: swelling
ratio slip P} lists the swelling ratios for increasing polynomial order, $P$,
for $B_{sl}=0.1$ and $B_{sl}=0.01$. The swelling ratios are converged to three
decimal places. Figure~\ref{subfig: newton swell slip fs u P} and 
\ref{subfig: newton swell slip fs v P} show that the velocity values
are converged for $P\geq 10$. The free surface spline for increasing slip
parameter is shown in Figure~\ref{fig:newt swell spline slip}. Increasing the
slip parameter yields a decrease in swelling.

\section{Conclusions}

In this article, we have demonstrated the capabilities of the high-order
spectral element method in the resolution of the stress singularity at the die
exit in the plane  Newtonian extrudate swell problem. We have shown that the
spectral method approximates the infinite pressure value with exponentially
increasing extreme values for increasing polynomial order. This high resolution
approximation of the steep stress profiles yields excellent predictions of the
swelling ratio. Our method predicts the same swelling ratio in comparison to low
order finite element methods with significantly fewer number of degrees of
freedom.

The only drawback of our high order method is the Gibbs oscillations, which
appear in the vicinity of the singularity for the pressure approximation. These
Gibbs oscillations are intrinsic to high order methods and they occur in the
approximation of discontinuities or steep profiles. However, we have
demonstrated that for the extrudate swell problem, the Gibbs oscillations stay
confined to one element next to the singularity and their amplitude decreases
significantly with increasing polynomial order. This small pollution in the
pressure profile is the price to pay in the high order method for the otherwise
excellent prediction of the steep pressure increase at the singularity.

We have given detailed results for a wide range of Reynolds numbers $0\leq \Rey
\leq 100$ in terms of swell ratios, exit pressure losses, free surface profiles
and velocity and pressure values. For the free surface profiles, we find three
extrudate swell regimes. The first is a reduction in swelling ($\Rey \leq 6$),
the second is a regime of a delayed swelling ($7 \leq \Rey \leq 10$) and the
third a contraction of the free liquid jet ($10 < \Rey \leq 100 $). With
increasing Reynolds number the maximum pressure values decrease and the Gibbs
oscillations decrease. We have then investigated the effect of slip along the
die wall. We have observed a reduction of the swelling for different slip
parameters $B_{sl}=\{0.01,0.1\}$ and have observed a drastic change in the
pressure profile which showed no occurrence of Gibbs oscillations.

\section{Acknowledgements}
The authors wish to warmly thank Prof. Evan Mitsoulis for generously providing
detailed data of his extrudate swell results.
The first author would like to thank the UK Engineering and Physical Sciences
Research Council for financial support.
\bibliographystyle{plainnat}


\end{document}